\documentclass{siamart190516}

\usepackage{graphics,amsmath,amsopn,amstext,amsfonts,color,fancybox,hyperref,bm,amssymb,cite}
\usepackage{booktabs,multirow,xspace}

\newcommand{\pppa}{\ensuremath{\texttt{PPPA}}\xspace}
\newcommand{\gurobi}{\ensuremath{\texttt{Barrier}}\xspace}

\newcommand{\revision}[1]{#1}

\newtheorem{remark}[theorem]{Remark}

\newcommand{\gap}{\vspace{0.1in}}

\newcommand{\epc}{\hspace{1pc}}
\newcommand{\thalf}{{\textstyle{\frac{1}{2}}}}

\newcommand{\onebld}{{\bf 1}}
\newcommand{\wt}{\widetilde}
\newcommand{\wh}{\widehat}

\title{Some Strongly Polynomially Solvable Convex Quadratic Programs with Bounded
Variables\footnote{Original March 2022; revised September 2022}}

\author{Jong-Shi Pang \thanks{Department of Industrial and
Systems Engineering, University of Southern California, Los Angeles, CA 90089
(\email{jongship@usc.edu; shaoning@usc.edu}).  This work was based on research partially
supported by the U.S.\ Air Force Office of Scientific Research under grant FA9550-18-1-0382.}
\and Shaoning Han \footnotemark[1]
}

\begin{document}

\maketitle



\begin{abstract}
\noindent This paper begins with the review of a class of strictly convex quadratic programs (QPs)
with bounded variables solvable by the parametric
principal pivoting algorithm with $\mbox{O}(n^3)$ strongly polynomial complexity,
where $n$ is the number of variables of the problem.
Extensions of this Hessian class are the main contributions of this paper, which
is motivated by a \revision{recent reference}
\cite{LiuFattahiGomezK21} wherein the efficient solution of a quadratic program
with a tridiagonal Hessian matrix in the quadratic objective
is needed for the construction of a polynomial-time algorithm for solving an associated
sparse variable selection problem.  With the tridiagonal structure, the
complexity of the QP algorithm reduces to $\mbox{O}(n^2)$.  Our strongly polynomiality
results extend
previous works of some strongly polynomially solvable linear complementarity
problems with a P-matrix \cite{PangChandrasekaran85}; special cases of the extended results
include weakly quasi-diagonally dominant problems in addition to the tridiagonal ones.
\end{abstract}

\section{Introduction} \label{sec:introduction}

\revision{The computational resolution of convex quadratic programs (QPs) has a long history; by resolution,
we mean either computing an optimal solution of the problem or determining that the problem is unbounded below,
assuming that it is feasible to begin with.  The early methods belong to several closely related families: simplex-type pivoting,
feasible direction, gradient projection, or active set.  In general, the computational
complexity of these methods have exponential worst-case behavior; see for instance the introduction in \cite{YeTse89}.
In contrast, the family of interior-point methods can solve convex QPs in polynomial time,
where the polynomiality is in terms of the number of variables and the binary bits to represent the size of the
problem data \cite{MonteiroAdler89,YeTse89}.}  This kind of polynomial complexity is different from the strongly
polynomial complexity wherein the polynomiality (in operation counts) is only in terms of the number of variables
and the number of constraints.  So far, QPs with the latter favorable solution
complexity are not well known in the literature, although special classes of
the related linear complementarity problems (LCPs) are known to be solvable by strongly
polynomially bounded algorithms.  The first such class of LCPs was discovered by
Chandraskaran \cite{Chandrasekaran70} for problems with a Z-matrix (i.e., with
nonpositive off-diagonal elements).  Subsequent extensions include
\cite{PangChandrasekaran85,AdlerCottlePang16}.  Most recently,
some quadratic sparse optimization problems with Stieltjes matrices (i.e., symmetric
positive definite Z-matrices) and bounded variables are shown to be strongly polynomially solvable
\cite{GomezHePang21,HeHanCuiGomezPang21}.  Till now, it is not known if there are bounded-variable
convex (but not strictly convex) QPs with non-Z Hessians that can be resolved in strongly polynomial time.

\gap

The present paper is partly motivated by the above void and also partly by
a \revision{recent reference} \cite{LiuFattahiGomezK21}.  The latter paper addresses the solution
of a special
class of quadratic sparse optimization problems with tridiagonal Hessian matrices
by a shortest-path algorithm with strongly polynomial complexity.
Nevertheless, the construction of the arc costs of such a network
formulation requires the efficient solution of quadratic programs in continuous
variables.  In the reference, the latter step is restricted to an unconstrained problem.
These two background sources lead us to the present study wherein we identify several
classes of convex QPs with bounded variables that can be resolved in
strongly polynomial time, where resolution means that either the successful computation
of an optimal solution or the demonstration that no such solution exists.

\section{Problem Formulation and the Positive Definite Case}

Consider the upper-and-lower-bounded quadratic program (QP) in the variable
$x \in \mathbb{R}^n$;
\begin{equation} \label{eq:QP}
\displaystyle{
\operatornamewithlimits{\mbox{\bf minimize}}_{0 \leq x \leq u}
} \ \theta(x) \, \triangleq \, q^{\top}x + \thalf \, x^{\top}Mx	
\end{equation}
where $0 < u_i \leq \infty$ for all $i \in [ n ] \triangleq \{ 1, \cdots, n \}$,
$q \in \mathbb{R}^n$ and $M \in \mathbb{S}^{\, n}$ is a symmetric matrix.
One way to solve (\ref{eq:QP}) is by solving the parametric problem:
\begin{equation} \label{eq:parametric QP}
\displaystyle{
\operatornamewithlimits{\mbox{\bf minimize}}_{0 \leq x \leq u}
} \ ( \, q + \tau \, p \, )^{\top}x + \thalf \, x^{\top}Mx,
\end{equation}
for a given vector $p$ so that $q + \tau_0 p \geq 0$ for some scalar $\tau_0 > 0$
and by tracing a solution path of (\ref{eq:parametric QP}) as a
function of the scalar $\tau$ via ``pivoting down'' $\tau$ until
$\tau = 0$ if possible, at which value a solution of the original problem (\ref{eq:QP})
is obtained.  In general, with $M$ being positive semidefinite, this solution strategy
can be accomplished by the classical parametric principal pivoting algorithm (PPPA)
for solving parametric linear complementarity problems (LCPs)
\cite[Algorithm~4.5.2]{CottlePangStone92}, thus for the convex QP (\ref{eq:QP}) in particular.
Although there is no formal result asserting the complexity of the strategy, the number of pivots
in the PPPA may be exponential in general when the vector $p$ is arbitrary.
Presented in \cite[Algorithm~4.8.2]{CottlePangStone92} the PPPA admits a simplified implementation
for an LCP with a P-matrix (that is not required to be symmetric).
\revision{By definition, a P-matrix is a real square matrix whose principal minors are all positive; there
are several characterizations of a P-matrix; see \cite[Theorem~3.3.4]{CottlePangStone92}.  We follow the standard
notation of submatrices and subvectors indexed by subsets of $[ n ]$; for instance, if $\alpha$ and $\beta$ are two
such index sets, then $M_{\alpha \beta}$ is the submatrix of $M$ with rows and columns indexed by elements in
$\alpha$ and $\beta$, respectively; similar definition applies to a subvector $q_{\alpha}$ of $q$.}
With the parametric vector $p$,
called an {\sl $n$-step vector}, satisfying the condition:
\begin{equation} \label{eq:n-step vector}
M_{\alpha \alpha}^{-1} p_{\alpha} \, \geq \, 0, \ \forall \, \alpha \, \subseteq \, [ n ], \epc \mbox{where} \epc
M_{\alpha \alpha}^{-1} \, \triangleq \, [ \, M_{\alpha \alpha} \, ]^{-1},
\end{equation}
it is shown in \cite[Theorem~4.8.3]{CottlePangStone92}, which
has its origin from \cite{PangChandrasekaran85}, that the
PPPA for an LCP with a P-matrix terminates in no more than $n$ pivots, thus justifying the terminology of an
$n$-step vector for $p$ and $n$-step scheme for the simplified algorithm.
Notice that such an $n$-step vector (if it exists) and the resulting algorithm are applicable
to all vectors $q$ in the linear term of the objective function of (\ref{eq:QP}) if $M$ is positive definite
(denoted $M \in \mathbb{S}_{++}^{\, n}$)
and the problem has no upper bounds.
An extension to problems with finite upper bounds is described in \cite[Algorithm~III]{GomezHePang21}
with the same linear-step convergence proved in \cite[Theorem~17]{GomezHePang21}. For convenience
of reference, we present the algorithm below.

\gap

\rule{6.4in}{0.01in}

\gap

{\bf Algorithm I: The $2n$-step PPPA for (\ref{eq:QP}) with $M$ positive definite}

\rule{6.4in}{0.01in}

\gap

\noindent {\bf Input.}  A matrix $M \in \mathbb{S}_{++}^{\, n}$ and a positive vector $p$ satisfying (\ref{eq:n-step vector}).

\gap

\noindent {\bf Initialization.}  Let $\left( \, \alpha, \, \beta, \, \gamma \, \right) =
\left( \, \emptyset, \, \{ 1, \cdots, n \}, \, \emptyset \, \right)$.  Also let $\tau_{\rm old} \triangleq \infty$.

\gap

\noindent {\bf General iteration.}  Solve a system of linear equations:
\begin{equation} \label{eq:barpq-basic}
M_{\alpha \alpha} \, ( \, \bar{q}_{\alpha},\bar{p}_{\alpha} \, ) \, = \,
( \, q_{\alpha} + M_{\alpha \gamma} u_{\gamma}, \, p_{\alpha} \, )
\end{equation}
for $( \, \bar{q}_{\alpha},\bar{p}_{\alpha} \, )$ and compute:
\begin{equation} \label{eq:barpq-nonbasic}
\begin{array}{lll}
( \, \bar{q}_{\beta},\bar{p}_{\beta} \, ) & \triangleq &
( \, q_{\beta} + M_{\beta \gamma} u_{\gamma}, \, p_{\beta} \, )
- M_{\beta \alpha} ( \, \bar{q}_{\alpha},\bar{p}_{\alpha} \, ) \\ [0.1in]
( \, \bar{q}_{\gamma},\bar{p}_{\gamma} \, ) & \triangleq &
( \, q_{\gamma} + M_{\gamma \gamma} u_{\gamma}, \, p_{\gamma} \, )
- M_{\gamma \alpha} ( \, \bar{q}_{\alpha},\bar{p}_{\alpha} \, ).
\end{array}
\end{equation}
Perform the ratio test (by convention, the maximum over an empty set is equal to $-\infty$):
\begin{equation} \label{eq:ratio test in pd}
\tau_{\rm new} \, \triangleq \, \max\left\{ \, \displaystyle{
\max_{i \in \beta}
} \, \left\{ \, -\displaystyle{
\frac{\bar{q}_i}{\bar{p}_i}
} \, \mid \, \bar{p}_i \, > \, 0 \right\}, \, \displaystyle{
\max_{i \in \alpha}
} \, \left\{ \, -\displaystyle{
\frac{u_i + \bar{q}_i}{\bar{p}_i}
} \, \mid \, \bar{p}_i \, > \, 0 \mbox{ and } u_i \, < \, \infty \, \right\}, \, 0 \, \right\}
\end{equation}
For $\tau \in \left[ \, \tau_{\rm new}, \, \tau_{\rm old} \, \right]$, let
$x_{\beta}(\tau) = 0$, $x_{\gamma}(\tau) = u_{\gamma}$, and $x_{\alpha}(\tau) = -\bar{q}_{\alpha} - \gamma \, \bar{p}_{\alpha}$.
Stop if $\tau_{\rm new} = 0$; in this case, $x(0)$ is the desired solution to (\ref{eq:QP}).  Otherwise
let $\bar{i} \in \alpha \cup \beta$ be a maximizing index of $\tau_{\rm new}$ (ties can be broken arbitrarily).  Update the index sets as follows:
\[
\left( \, \alpha, \, \beta, \, \gamma \, \right)_{\rm new} \, = \, \left\{ \begin{array}{ll}
\left( \, \alpha \, \cup \, \{ \, \bar{i} \, \}, \ \beta \, \setminus \, \{ \, \bar{i} \, \}, \ \gamma  \, \right) & \mbox{if $i \in \beta$} \\ [0.1in]
\left( \, \alpha \, \setminus \, \{ \bar{i} \}, \ \beta, \ \gamma \, \cup \, \{ \, \bar{i} \, \} \, \right) & \mbox{if $i \in \alpha$}.
\end{array} \right.
\]
Set $\tau_{\rm old} \leftarrow \tau_{\rm new}$ and return to the beginning of the general iteration.  \hfill $\Box$

\gap

\rule{6.4in}{0.01in}

\gap

We summarize the $2n$-step termination of the above algorithm; \revision{see the cited theorem in \cite{GomezHePang21}}.
A noteworthy point about this result is that there is
no nondegeneracy assumption needed.  This remark applies to all the extended strongly polynomial results in the rest of
the paper.

\begin{proposition} \label{pr:pd strongly poly} \rm
Let $M \in \mathbb{S}^{\, n}$ be positive definite.  If there exists a vector
$p \in \mathbb{R}_{++}^n$ satisfying (\ref{eq:n-step vector}), then
the unique solution of (\ref{eq:QP}) can be obtained by the PPPA in no more than $2n$ iterations.
\hfill $\Box$
%
\end{proposition}

\begin{remark} \rm
With $M$ being symmetric positive definite, there are some interesting necessary
and sufficient conditions for the existence of an $n$-step vector $p$; see
\cite[Corollary~4.8.11 and Theorem~3.11.18]{CottlePangStone92}.  In particular,
by the cited
Corollary, if such a vector $p$ exists, then the unique solution the QP (\ref{eq:QP})
with all $u_i = \infty$ has a \revision{geometric} characterization \revision{as a ``least element'' of
the feasible region of the Karush-Kuhn-Tucker conditions of (\ref{eq:QP}),
which is the set: $\{ x \geq 0 \mid q + Mx \geq 0 \}$,}
least in a certain cone ordering of $\mathbb{R}^n$.  \hfill $\Box$
\end{remark}

Taking advantage of the single-row and
(same) column changes in the working principal submatrix $M_{\alpha \alpha}$,
the two computations in (\ref{eq:barpq-basic}) and (\ref{eq:barpq-nonbasic}) can be
carried out in $\mbox{O}(n^2)$ (instead of the general
$\mbox{O}(n^3)$) algebraic operations; this is similar to the standard linear-algebraic
based implementation of the revised simplex method by basis updates; see the computational remark below.
Consequently, with the $2n$-step termination, the  overall parametric scheme computes
the unique solution of (\ref{eq:QP}) with $\mbox{O}(n^3)$ complexity when $M$ is positive
definite and an $n$-step vector exists.  This complexity can be further reduced to
O$(n^2)$ when $M$ is tridiagonal since the linear algebraic calculations can be carried out
in $\mbox{O}(n)$ time.

\gap

The primary goal of this paper is to extend the above algorithm and result to several subclasses of the class
$\mathbb{S}_+^{\, n}$ of symmetric positive semidefinite matrices.
%
%
Before leaving the positive definite case, we highlight two noteworthy points
about the positive semidefinite case: (i) the problem (\ref{eq:QP}) may have
an unbounded objective value when some upper bounds are infinite; and
(ii) computationally, certain principal submatrices $M_{\alpha \alpha}$ may be
singular, thus the condition
(\ref{eq:n-step vector}) immediately becomes invalid and requires modification.

\gap

\noindent \revision{
{\bf A computational remark (cf.\ \cite[Section~4.10]{CottlePangStone92}).}
Consider two systems of linear equations, where the index $\bar{i} \not\in \alpha$,
\[
M_{\alpha \alpha} \, \bar{p}_{\alpha} \, = \, p_{\alpha} \epc \mbox{and} \epc
\left[ \begin{array}{cc}
M_{\alpha \alpha} & M_{\alpha \bar{i}} \\ [5pt]
M_{\bar{i} \alpha} & m_{\bar{i} \bar{i}}
\end{array} \right] \, \left( \begin{array}{c}
( \bar{p}_{\alpha} )_{\rm new} \\ [5pt]
\bar{p}_{\bar{i}}
\end{array} \right) \, = \, \left( \begin{array}{c}
p_{\alpha} \\ [5pt]
p_{\bar{i}}
\end{array} \right)
\]
The relations between the solutions of these two systems are as follows:
\[ \begin{array}{cll}
( \bar{p}_{\alpha} )_{\rm new} & = & \bar{p}_{\alpha} - M_{\alpha \alpha}^{-1} M_{\alpha \bar{i}} \, \bar{p}_{\bar{i}} \\ [5pt]
\bar{p}_{\bar{i}} & = & \left( \, m_{\bar{i} \bar{i}} - M_{\bar{i} \alpha} M_{\alpha \alpha}^{-1} M_{\alpha \bar{i}} \, \right)^{-1}
\left( \, p_{\bar{i}} - M_{\bar{i} \alpha} \bar{p}_{\alpha} \, \right).
\end{array} \]
Suppose that $\bar{p}_{\alpha}$ is known and the index $\bar{i}$ is added to the current index set $\alpha$.  Then with a symmetric $M$,
provided that $M_{\alpha \alpha}^{-1}$ is available, the vector $M_{\alpha \alpha}^{-1} M_{\alpha \bar{i}}$, and thus the updated vector
$( \bar{p}_{\alpha} )_{\rm new}$ and the new component $\bar{p}_{\bar{i}}$, can be readily computed with some vector-matrix multiplications.
The total computational effort is $\mbox{O}(| \alpha |^2)$ arithmetic operations.  Moreover, we have
\[
\left[ \begin{array}{cc}
M_{\alpha \alpha} & M_{\alpha \bar{i}} \\ [0.1in]
M_{\bar{i} \alpha} & m_{\bar{i} \bar{i}}
\end{array} \right]^{-1} \, = \, \left[ \begin{array}{cc}
M_{\alpha \alpha}^{-1} + \displaystyle{
\frac{M_{\alpha \alpha}^{-1} M_{\alpha \bar{i}} M_{\bar{i} \alpha} M_{\alpha \alpha}^{-1}}{m_{\bar{i} \bar{i}} -
M_{\bar{i} \alpha} M_{\alpha \alpha}^{-1}M_{\alpha \bar{i}}}
} & \displaystyle{
\frac{-M_{\alpha \alpha}^{-1}M_{\alpha \bar{i}}}{m_{\bar{i} \bar{i}} - M_{\bar{i} \alpha} M_{\alpha \alpha}^{-1}M_{\alpha \bar{i}}}
} \\ [0.2in]
\displaystyle{
\frac{- M_{\bar{i} \alpha} M_{\alpha \alpha}^{-1}}{m_{\bar{i} \bar{i}} - M_{\bar{i} \alpha} M_{\alpha \alpha}^{-1}M_{\alpha \bar{i}}}
} & \displaystyle{
\frac{1}{m_{\bar{i} \bar{i}} - M_{\bar{i} \alpha} M_{\alpha \alpha}^{-1}M_{\alpha \bar{i}}}
}
\end{array} \right],
\]
which shows that the inverse of the left-hand bordered matrix can similarly be computed with an available $M_{\alpha \alpha}^{-1}$ in
$\mbox{O}(| \alpha |^2)$ arithmetic operations.
The updating in the case where an index $\bar{i}$ is dropped from the current index set $\alpha$ can be accomplished by reversing the
above calculations.  Therefore, provided that we maintain the inverse $M_{\alpha \alpha}^{-1}$ throughout the calculations,
the linear-algebraic steps in each iteration of Algorithm~I can be accomplished in
$\mbox{O}(n^2)$ operations at the most.  When the given matrix $M$ is sparse and instead of manipulating the inverse $M_{\alpha \alpha}^{-1}$,
advanced matrix factorization techniques used in the practical implementation of the simplex method can be applied to carry out the
vector-matrix updates without solving the systems of linear equations from scratch at each iteration.   \hfill $\Box$
}

\section{Comparison Matrices and Schur Complements} \label{sec:comparison and Schur}

For a given matrix $M \in \mathbb{R}^{n \times n}$,
its {\sl comparison matrix} denoted $\overline{M}$, has entries given by
\begin{equation} \label{eq:entries comparison}
\overline{M}_{ij} \, \triangleq \, \left\{ \begin{array}{cl}
m_{ii} & \mbox{if $i = j$} \\ [5pt]
- | \, m_{ij} \, | & \mbox{if $i \neq j$}.
\end{array} \right.
\end{equation}
Notice that we employ the original diagonal entries $m_{ii}$ of $M$ to define the same
entries of $\overline{M}$.  This allows us to omit the phrase ``with nonnegative diagonal
entries'' throughout the rest of the paper.  Obviously, $\overline{M}$ is a Z-matrix,
i.e., its off-diagonal entries are all nonpositive.  Schur complements \cite{Cottle74} play an important role in this work.  Specifically,
given a nonsingular principal submatrix of $M$, say
$M_{\alpha \alpha}$ where $\alpha \subseteq [ n ]$ with complement $\bar{\alpha}$,
a principal pivot on $M_{\alpha \alpha}$ transforms the linear system:
\[ \begin{array}{lll}
w_{\alpha} & = & r_{\alpha} + M_{\alpha \alpha} x_{\alpha} + M_{\alpha \bar{\alpha}}	
x_{\bar{\alpha}} \\ [0.1in]
w_{\bar{\alpha}} & = & r_{\bar{\alpha}} + M_{\bar{\alpha} \alpha} x_{\alpha} +
M_{\bar{\alpha} \bar{\alpha}} x_{\bar{\alpha}}
\end{array} \]
into the equivalent system:
\begin{equation} \label{eq:principal pivot transform in full}
\begin{array}{lll}
x_{\alpha} & = & \epc -M_{\alpha \alpha}^{-1} r_{\alpha} \hspace{0.35in} + \epc
M_{\alpha \alpha}^{-1} w_{\alpha} \hspace{0.1in} - \
M_{\alpha \alpha}^{-1} M_{\alpha \bar{\alpha}}
x_{\bar{\alpha}} \\ [0.1in]
w_{\bar{\alpha}} & = & r_{\bar{\alpha}} - M_{\bar{\alpha} \alpha} M_{\alpha \alpha}^{-1}
r_{\alpha} + M_{\bar{\alpha} \alpha} M_{\alpha \alpha}^{-1} w_{\alpha} +
( M/M_{\alpha \alpha} ) x_{\bar{\alpha}}
\end{array}
\end{equation}
where $( M/M_{\alpha \alpha} )$ is the {\sl Schur complement} of $M_{\alpha \alpha}$ in
$M$ given by:
\[
\left( M/M_{\alpha \alpha} \, \right) \, \triangleq \, M_{\bar{\alpha} \bar{\alpha}} -
M_{\bar{\alpha} \alpha} M_{\alpha \alpha}^{-1} M_{\alpha \bar{\alpha}}.
\]
The vector and matrix in
(\ref{eq:principal pivot transform in full}):
\begin{equation} \label{df:principal pivot transform}
\left( \begin{array}{c}
-M_{\alpha \alpha}^{-1} r_{\alpha} \\ [0.1in]
r_{\bar{\alpha}} - M_{\bar{\alpha} \alpha} M_{\alpha \alpha}^{-1}
r_{\alpha}
\end{array} \right) \epc \mbox{and} \epc
\left[ \begin{array}{cc}
M_{\alpha \alpha}^{-1} & -M_{\alpha \alpha}^{-1} M_{\alpha \bar{\alpha}} \\ [0.1in]
M_{\bar{\alpha} \alpha} M_{\alpha \alpha}^{-1} & ( M/M_{\alpha \alpha} )
\end{array} \right]
\end{equation}
are called the {\sl principal pivot transforms} of $r$ and $M$ with reference to
$M_{\alpha \alpha}$, respectively.  Note that if $M$ is symmetric, then so is $( M/M_{\alpha \alpha} )$.
The well-known Schur determinantal formula states that
\[
\det(M) \, = \, \det( M_{\alpha \alpha} ) \, \det( M/M_{\alpha \alpha} ).
\]
This formula yields in particular the following result asserting the positive semidefiniteness of $M$
in terms of that of its Schur complement.

\begin{proposition} \label{pr:psd in terms of Schur} \rm
Let $M \in \mathbb{S}^{\, n}$.  The following two statements hold.

\gap

(a) If $M$ is positive semidefinite and has a positive definite principal submatrix $M_{\alpha \alpha}$,
then the Schur complement $( M/M_{\alpha \alpha} )$ is positive semidefinite.

\gap

(b) Conversely, if $M$ has a positive diagonal element $m_{ii} > 0$ and the Schur complement $( M/m_{ii} )$ is
positive semidefinite, then $M$ is positive semidefinite.
\end{proposition}

\begin{proof}  ``Only if''.  Let $\beta \triangleq [ n ] \setminus \alpha$.  To show that $( M/M_{\alpha \alpha} )$ is positive semidefinite,
it suffices to show that $\det\left( ( M/M_{\alpha \alpha} )_{\beta^{\, \prime} \beta^{\, \prime}} \right)$ is nonnegative for all subsets
$\beta^{\, \prime}$ of $\beta$.  This holds easily because $( M/M_{\alpha \alpha} )_{\beta^{\, \prime} \beta^{\, \prime}}$
is the Schur complement of $M_{\alpha \alpha}$ in $\left[ \begin{array}{cc}
M_{\alpha \alpha} & M_{\alpha \beta^{\, \prime}} \\ [5pt]
M_{\beta^{\, \prime} \alpha} & M_{\beta^{\, \prime} \beta^{\, \prime}}
\end{array} \right]$.  Thus the nonnegativity of $\det\left( ( M/M_{\alpha \alpha} )_{\beta^{\, \prime} \beta^{\, \prime}} \right)$
follows readily from the Schur determinantal formula applied to the latter $( \alpha \cup \beta^{\, \prime} )$-principal submatrix.

\gap

Conversely, suppose that $m_{ii} > 0$ and the Schur complement $( M/m_{ii} )$ is positive semidefinite.  Let $\beta = [ n ] \setminus \{ i \}$.
Then
\[
M_{\beta \beta} \, = \, ( M/m_{ii} ) + \displaystyle{
\frac{M_{\beta i} M_{i \beta}}{m_{ii}},
} \]
where the right-hand side is the sum of a positive semidefinite matrix and a symmetric rank-one matrix.  Hence $M_{\beta \beta}$ is positive semidefinite.
To complete the proof, it suffices to show that for any $\beta^{\, \prime} \subseteq \beta$,
\[
\det\left( \left[ \begin{array}{cc}
m_{ii} & M_{i \beta^{\, \prime}} \\ [5pt]
M_{\beta^{\, \prime} i}  & M_{\beta^{\, \prime} \beta^{\, \prime}}
\end{array} \right] \right) \, \geq \, 0.
\]
The left-hand determinant is equal to the ratio:
\[
\displaystyle{
\frac{1}{m_{ii}}
} \, \det\left( M_{\beta^{\, \prime} \beta^{\, \prime}} - \displaystyle{
\frac{M_{\beta^{\, \prime} i} M_{i \beta^{\, \prime}}}{m_{ii}}
} \right) \, = \, \displaystyle{
\frac{1}{m_{ii}}
} \, \det\left( ( M/m_{ii} )_{\beta^{\, \prime} \beta^{\, \prime}} \right) \, \geq \, 0
\]
because $( M/m_{ii} )_{\beta^{\, \prime} \beta^{\, \prime}}$ is a principal submatrix of the positive semidefinite $( M/m_{ii} )$.
\end{proof}

\gap

It would be useful to summarize some well-known necessary and sufficient
conditions for a symmetric Z-matrix to be positive semidefinite.
\revision{In particular, the vector $d$ satisfying condition (c) can be shown to be an $n$-step vector
for a Stieltjes matrix; this condition provides the basis for the generalization of such a vector to the class of matrices
$\overline{\mathbb{S}}^{\, n}_+$ to be introduced in the next section.}

\begin{lemma} \label{lm:symmetric Z psd} \rm
Let $A \in \mathbb{S}^{\, n}$ be a Z-matrix.  The following statements are
equivalent.

\gap

(a) $A$ is positive semidefinite;

\gap

(b) $A$ is a P$_0$-matrix, that is, all its principal minors are nonnegative;

\gap

(c) there exists a vector $d > 0$ such that $Ad \geq 0$. 	
\end{lemma}

\begin{proof}  (a) $\Leftrightarrow$ (b).  That the P$_0$-property is equivalent to
positive semidefiniteness of any symmetric matrix is well known.

\gap

(b) $\Leftrightarrow$ (c).  This follows from \cite[Theorem~3.4]{PooleBoullion74}
if $A$ is irreducible.  If $A$ is reducible, then it must be block diagonal.   An
inductive argument then easily completes the proof of the equivalence of (b) and (c).
\end{proof}

Property (c) asserts the {\sl weakly quasi-diagonally dominance} of the Z-matrix
$A$; indeed, the condition $Ad \geq 0$ with $A$ being a Z-matrix and $d > 0$
is equivalent to
\[
a_{ii} \, d_i \, \geq \, \displaystyle{
\sum_{j \neq i}
} \, | \, a_{ij} \, | \, d_j \epc \forall \, i \, \in \, [ n ].
\]
%
%
It is known that if $\overline{M} \in \mathbb{S}^{\, n}$ is positive definite, (thus is a
Stieltjes matrix), then for any vector $d \in \mathbb{R}^n_{++}$ for which
$\overline{M}d > 0$ (such a vector must
exist), the vector $p = \thalf \, ( \, M + \overline{M} )d$ is positive and satisfies
(\ref{eq:n-step vector}); see \cite[Corollary~1]{PangChandrasekaran85}.   This result
is the basis for extension to the positive semidefinite case.

\section{The Matrix Class $\overline{\mathbb{S}}^{\, n}_+$}

We denote the class of matrices $M \in \mathbb{S}^{\, n}$ for which $\overline{M}$ is
positive semidefinite by $\overline{\mathbb{S}}^{\, n}_+$.  With the definition
(\ref{eq:entries comparison}) of the diagonal entries of $\overline{M}$, it follows that
a matrix $M \in \overline{\mathbb{S}}^{\, n}_+$ must have nonnegative diagonal entries.
Part (a) of the result below
establishes the inclusion $\overline{\mathbb{S}}^{\, n}_+ \subseteq
\mathbb{S}^{\, n}_+$; parts (b) and (c) show that the positive semidefiniteness
of the comparison matrix $\overline{M}$ of $M$ is inherited by the comparison matrices
of all principal submatrices of $M$ and those of its
Schur complements.  \revision{These two parts are instrumental for a reduction process to be
discussed in Section~\ref{sec:extended problem classes}}.
Part (d) identifies a special vector $p$ that arises from
the positive semidefinitenss of $\overline{M}$.
This vector in the special case of a matrix $M \in \overline{\mathbb{S}}^{\, n}_+$ provides the background
for the discussion in Section~\ref{sec:PPPA with p} for a general positive semidefinite matrix.

\begin{lemma} \label{lm:schur psd} \rm
Let $M \in \overline{\mathbb{S}}^{\, n}_+$ be given.  The following four statements hold:

\gap

(a) $M$ is itself positive semidefinite; thus, $\overline{\mathbb{S}}^{\, n}_+ \subseteq
\mathbb{S}^{\, n}_+$;

\gap

(b) the comparison matrices of all principal submatrices of $M$ are positive
semidefinite;

\gap

(c) if $M_{\alpha \alpha}$ is a nonsingular principal submatrix of $M$ with
$\alpha \subseteq [ n ]$, then the comparison matrix of the Schur complement
$\left( M/M_{\alpha \alpha} \, \right)$ is positive semidefinite;

\gap

(d) for any vector $d \in \mathbb{R}^n_{++}$ for which
$\overline{M}d \geq 0$ (such a vector must
exist by part (c) of Lemma~\ref{lm:symmetric Z psd}), the vector
$p = \thalf \, ( \, M + \overline{M} )d$ is nonnegative and satisfies the following
two conditions:

\gap

$\bullet $ $p_i = 0$ only if the off-diagonal entries of the $i$th row of $M$
are all nonpositive (thus, $p > 0$ if every row of $M$ has at least one positive off-diagonal element);
moreover,

\gap

$\bullet $ for any nonempty index set $\alpha \subset [ n ]$ and index
$i \not\in \alpha$,
\begin{equation} \label{eq:psd n-step vector}
[ M_{\alpha \alpha} ]^{-1} \ \mbox{ exists } \, \Rightarrow \,
[ \, M_{\alpha \alpha} \, ]^{-1} p_{\alpha} \, \geq \, 0 \ \mbox{ and } \
p_i - M_{i \alpha} [ \, M_{\alpha \alpha} \, ]^{-1} p_{\alpha} \, \geq \, 0.
\end{equation}
\end{lemma}

\begin{proof} Since $M$ has nonnegative diagonal entries, its positive semidefiniteness
is immediate from the inequality:
\[
x^{\top}Mx \, \geq \, | \, x \, |^{\top} \, \overline{M} \, | \, x \, |,
\epc \forall \, x \, \in \, \mathbb{R}^n,
\]	
where $| \, x \, |$ is the vector whose components are the absolute values of those
of $x$.  This shows that $M$ is positive semidefinite.  Statement (b)
is obvious.  To prove (c), let
$M_{\alpha \alpha}$ be a nonsingular principal submatrix of $M$, where
$\alpha \subseteq [ n ]$.  It follows that $M_{\alpha \alpha}$ is positive definite.
To show that
$\overline{ \left( M/M_{\alpha \alpha} \right)}$ is positive semidefinite,
consider obtaining the Schur complement $(  M/M_{\alpha \alpha} )$ by Gaussian pivoting, i.e.,
successively pivoting on diagonal elements.  To be precise,
let $\alpha \triangleq \{ 1, \cdots, K \}$.  Let
$M^1 \triangleq \left( M/m_{11} \right)$ and $m^{\, 1}_1$ be the first diagonal
element of $M^{\, 1}$.  Inductively, define
\[ \left\{ \begin{array}{lll}
M^k & \triangleq & \left( M^{k-1}/m^{\, k-1}_1 \right) \\ [5pt]
m^{\, k}_1 & \triangleq & \mbox{first diagonal element of $M^k$}
\end{array} \right\} \epc k = 2, \cdots, K-1.
\]
With each $M^k$ being positive definite for $k = 1, \cdots, K-1$, we can show that the Schur complement
$\left( M/M_{\alpha \alpha} \, \right) = \left( \, M^{K-1}/m^{\, K-1}_1 \, \right)$.  With this
identification of $\left( M/M_{\alpha \alpha} \, \right)$, to establish the
positive semidefinitenss of the comparison matrix of
$\left( M/M_{\alpha \alpha} \, \right)$, it suffices to consider the case when
$\alpha$ is a singleton, say $\alpha = \{ 1 \}$ and $m_{11} > 0$.  We need to show the inequality in
\[
y^{\top} \, \overline{\left( M/M_{\alpha \alpha} \, \right)} \, y \, = \, \displaystyle{
\sum_{i=2}^n
} \, y_i^2 \, \left( \, \underbrace{m_{ii} - \displaystyle{
\frac{m_{i1}^2}{m_{11}}
}}_{\mbox{nonnegative}} \, \right) - \displaystyle{
\sum_{\substack{i,j=2 \\ i \neq j}}^n
} \, y_i \, y_j \, \left| \, m_{ij} - \displaystyle{
\frac{m_{i1} \, m_{1j}}{m_{11}}
} \, \right| \, \geq \, 0, \epc \forall \, y = ( y_i )_{i=2}^n.
\]
Since
\[
\left| \, m_{ij} - \displaystyle{
\frac{m_{i1} \, m_{1j}}{m_{11}}
} \, \right| \, \leq \, | \, m_{ij} \, | + \displaystyle{
\frac{| \, m_{i1} \, | \, | \, m_{1j} \, |}{m_{11}}
} \, ,
\]
it follows that
\[ \begin{array}{lll}
y^{\top} \, \overline{\left( M/M_{\alpha \alpha} \, \right)} \, y & \geq &
\displaystyle{
\sum_{i=2}^n
} \, y_i^2 \, \left( \, m_{ii} - \displaystyle{
\frac{m_{i1}^2}{m_{11}}
} \, \right) - \displaystyle{
\sum_{\substack{i,j=2 \\ i \neq j}}^n
} \, | \, y_i \, | \, | \,  y_j \, | \, \left( \, | \, m_{ij} \, | + \displaystyle{
\frac{| \, m_{i1} \, | \, | \, m_{1j} \, |}{m_{11}}
} \, \right) \\ [0.2in]
& = & | \, y \, |^{\top} \left( \, \overline{M}/m_{11} \, \right) | \, y \, |
\, \geq \, 0,
\end{array} \]
where the last inequality holds because $\left( \, \overline{M}/m_{11} \, \right)$
is the Schur complement of $m_{11}$ in the positive semidefinite matrix $\overline{M}$.

\gap

For part (d), we first show that the vector $p$ is nonnegative.  This follows easily
because
\[
p_i \, = \, m_{ii}d_i + \thalf \, \displaystyle{
\sum_{j \neq i}
} \, \left[ \, m_{ij} - | \, m_{ij} \, | \, \right] d_j \, \geq \,
m_{ii}d_i - \displaystyle{
\sum_{j \neq i}
} \, | \, m_{ij} \, | \, d_j \, = \, ( \, \overline{M} d \, )_i \, \geq \, 0.
\]
Moreover, if $p_i = 0$, then we must have $m_{ij} \leq 0$ for all $j \neq i$.  To prove
(\ref{eq:psd n-step vector}), we note that for every scalar $\varepsilon > 0$,
the comparison matrix of the perturbed
matrix $M^{\, \varepsilon} \triangleq M + \varepsilon \, \mathbb{I}$,
where $\mathbb{I}$ is the identity matrix, is positive definite.  Hence,
it follows that the vector
$p^{\, \varepsilon} \triangleq \displaystyle{
\frac{( \, M + \overline{M} \, ) d}{2}
} + \varepsilon \, d$ satisfies
\[
\left[ \, M_{\alpha \alpha} + \varepsilon \, \mathbb{I} \, \right]^{-1}
p_{\alpha}^{\, \varepsilon} \, > \, 0, \epc \forall \, \alpha \, \subseteq \, [ n ].
\]
Take a sequence of positive scalars $\{ \varepsilon_k \}$ converging to zero.
The corresponding sequence
$\{ \, p^{\, \varepsilon_k} \, \}$ converges to the vector $p$.
If $M_{\alpha \alpha}$ is invertible, then it follows that
$[ \, M_{\alpha \alpha} \, ]^{-1}p_{\alpha} \geq 0$.  Finally, to prove the nonnegativity
of $p_i - M_{i \alpha} [ \, M_{\alpha \alpha} \, ]^{-1}p_{\alpha}$ for $i \not\in \alpha$, let
\[
\left( \begin{array}{c}
\bar{p}_{\alpha}^{\, \varepsilon} \\ [5pt]
\bar{p}_i^{\, \varepsilon}	
\end{array} \right) \, \triangleq \,\left[ \begin{array}{ccc}
M_{\alpha \alpha} + \varepsilon \, \mathbb{I} & M_{\alpha i} \\ [5pt]
M_{i \alpha} & m_{ii} + \varepsilon
\end{array} \right]^{-1} \left( \begin{array}{c}
p_{\alpha}^{\, \varepsilon} \\ [5pt]
p_i^{\, \varepsilon}	
\end{array} \right) \, > \, 0,
\]
from which we obtain the equation:
\[
p_i^{\, \varepsilon} - M_{i \alpha} \left[ \, M_{\alpha \alpha} +
\varepsilon \, \mathbb{I} \, \right]^{-1} p^{\, \varepsilon}_{\alpha} \, = \, \left( \,
\underbrace{m_{ii} + \varepsilon - M_{i \alpha} \left[ \, M_{\alpha \alpha} +
\varepsilon \, \mathbb{I} \, \right]^{-1} M_{\alpha i}}_{\mbox{posiitve}} \, \right) \,
\bar{p}^{\, \varepsilon}_i.
\]
Therefore, the left-hand side is nonnegative for all $\varepsilon > 0$.  Letting
$\varepsilon \downarrow 0$ completes the proof.
\end{proof}

The next result identifies some matrices in $\overline{\mathbb{S}}^{\, n}_+$ and
has three parts.  \revision{With relevance to the discussion
in Subsection~\ref{subsec:reducibility and tridiagonality} about
(\ref{eq:QP}) with a tridiagonal matrix $M$, which is a class of QPs in
\cite{LiuFattahiGomezK21} that has motivated our study.
Although parts (b) and (c) do not have direct implications to the
subsequent development, they are included here to demonstrate the breadth and some basic properties
of this class of matrices}.  In particular, part (b) shows that the sum of two matrices in
the class $\overline{\mathbb{S}}^n_+$ remains a matrix in this class; part (c)
concerns a special block structured matrix one of whose diagonal blocks is tridiagonal.

\begin{proposition} \label{pr:some further subclasses} \rm
The following three statements hold for a matrix $M \in \mathbb{S}^{\, n}$.

\gap

(a) If $M$ is tridiagonal, then $\overline{M}$ is
positive semidefinite if and only if $M$ is.

\gap

(b)  Let $N \in \mathbb{S}^{\, n}$. 
If both $\overline{M}$ and $\overline{N}$ are positive semidefinite,
then so is $\overline{M + N}$.

\gap

(c)  Suppose that $M \in \mathbb{S}_+^{\, n}$ and for some $k \in [ n ]$,
$M_{\alpha \alpha} \in \overline{\mathbb{S}}_+^{\, k-1}$ and
$M_{\gamma \gamma}$ is tridiagonal, where $\alpha = \{ 1, \cdots, k-1 \}$ and $\gamma = \{ k, \cdots, n \}$,
moreover $M_{\alpha \gamma}$ is the zero matrix except for the entry $m_{k-1 \, k}$ which is nonzero.
Then $M \in \overline{\mathbb{S}}_+^{\, n}$.
\end{proposition}

\begin{proof}
To prove statement (a), it suffices to show that if $M$ is tridiagonal positive
semidefinite, then $\overline{M}$ is positive semidefinite.  If $m_{11} = 0$, then
the entire first row and same column of $M$ is equal to zero, by its
positive semidefiniteness; an inductive argument then yields the
positive semidefiniteness of $\overline{M}$ immediately.  So we may assume that
$m_{11} > 0$ and consider a principal pivot on $m_{11}$.  Since $M$ is
tridiagonal and symmetric, we have
\[
\left( \overline{M}/m_{11} \right) \, = \, \left\{ \begin{array}{l}
\mbox{matrix $\overline{M}$ with the first row and column deleted} \\ [5pt]
\mbox{and the $m_{22}$ element changed to $\bar{m}_{22} \triangleq m_{22} -
\displaystyle{
\frac{m_{12}^2}{m_{11}}
}$}
\end{array} \right\} \, = \, \overline{\left( M/m_{11} \right)},
\]
where $\left( \overline{M}/m_{11} \right)$ is the Schur complements of $m_{11}$ in
$\overline{M}$ and $\overline{\left( M/m_{11} \right)}$ is the comparison matrix of
the Schur complement of $m_{ii}$ in $M$.  Note that the modified diagonal entry
$\bar{m}_{22}$ is
nonnegative by the positive semidefiniteness of $M$.  Since $( M/m_{11} )$ is a
tridiagonal positive semidefinite matrix of order $n-1$, an induction
hypothesis yields that its comparison matrix is positive semidefinite; thus so is
$\left( \overline{M}/m_{11} \right)$.  By part (b) of Proposition~\ref{pr:psd in terms of Schur}, the
positive semidefiniteness of $\overline{M}$ follows.

\gap

To prove statement (b), it suffices to note that for all pairs $(i,j)$,
\[
-| \, m_{ij} + n_{ij} \, | \, \geq \, -| \, m_{ij} \, | - | \, n_{ij} \, |.
\]
Hence, $\overline{M + N} \geq \overline{M} + \overline{N}$.  Consequently,
\begin{equation} \label{eq:key inequality}
x^{\top} \left( \, \overline{M + N} \, \right) x \, \geq \,
| \, x \, |^{\top} \left( \, \overline{M + N} \, \right) | \, x \, | \, \geq \,
| \, x \, |^{\top} \, \overline{M} \, | \, x \, | + | \, x \, |^{\top} \, \overline{N} \,
| \, x \, |, \epc \forall \ x.
\end{equation}
Hence $\overline{M + N}$ is clearly positive semidefinite if both
$\overline{M}$ and $\overline{N}$ are.

\gap

For part (c), we may assume that $m_{kk} > 0$.
The augmented principal submatrix $M_{\alpha_1 \alpha_1}$,
where $\alpha_1 = \alpha \cup \{ k \}$, has the structure
\[
M_{\alpha_1 \alpha_1} \, = \, \left[ \begin{array}{cc}
M_{\alpha \alpha} & a \\ [0.1in]
a^{\top} & m_{kk}
\end{array} \right]
\]
where $a \in \mathbb{R}^{k-1}$ is the zero vector except for its last element which is equal to $m_{k-1 \, k}$.
Since the last row and same column of $M_{\alpha_1 \alpha_1}$ has the same structure of a tridiagonal matrix,
an argument similar to part (a) shows that $M_{\alpha_1 \alpha_1}$ belongs to $\overline{\mathbb{S}}_+^k$.
Next we consider $M_{\alpha_2 \alpha_2}$, where $\alpha_2 = \alpha_1 \cup \{ k + 1 \}$.  We have
\[
M_{\alpha_2 \alpha_2} \, = \, \left[ \begin{array}{cc}
M_{\alpha_1 \alpha_1} & b \\ [0.1in]
b^{\top} & m_{k+1 \, k+1}
\end{array} \right],
\]
where $b \in \mathbb{R}^{k}$ is the zero vector except for its last element which is equal to $m_{k \, k+1}$.
A similar argument shows that $M_{\alpha_2 \alpha_2} \in \overline{\mathbb{S}}_+^{k+1}$.  Continuing this way,
we can finally establish that $M \in \overline{\mathbb{S}}_+^{\, n}$.
\end{proof}

\begin{remark} \rm
In general, $\overline{M + N} \neq \overline{M} + \overline{N}$ and the inequality
(\ref{eq:key inequality}) is the key in the above proof of statement (b).  \hfill $\Box$
\end{remark}

\subsection{The irreducible subclass} \label{subsec:irreducible}

\revision{Obviously, if $M$ is a block diagonal matrix, then the problem (\ref{eq:QP}) decomposes into
smaller-sized subproblems each of the same kind.  When $M$ is symmetric, block diagonality is equivalent to reducibility,
a well-known matrix-theoretic property whose definition we recall.}
A matrix $A \in \mathbb{R}^n$ is irreducible if there does not exist
a permutation matrix $P$ such that
\[
PAP^{\top} \, = \, \left[ \begin{array}{ll}
B & 0 \\ [5pt]
C & D
\end{array} \right],
\]
where $B$ and $D$ are square matrices; a reducible matrix is one that is not irreducible.  It follows that a matrix
$M \in \mathbb{S}^{\, n}$ is \revision{reducible} if and only if it is block diagonal.  In particular,
a tridiagonal $M \in \revision{\mathbb{S}^{\, n}_+}$ is irreducible if and only if its diagonal and super-diagonal
(and thus the corresponding lower-diagonal) have no zero entries.
It is known \cite[6.2.2]{Ortega90} that a matrix $A$ is irreducible if and only if
for any two distinct indices $i \neq j$, there exists a chain of nonzero entries of $A$ of the form:
$\{ \, a_{i_1 i}, \cdots, a_{i_m j} \}$.  In general, it is a classic fact that checking if an $n \times n$ matrix
is reducible or not
can be accomplished by a graph-based algorithm with O$(n^2)$ complexity.
In what follows, we establish some basic properties of an irreducible
matrix $A \in \overline{\mathbb{S}}_+^{\, n}$.  \revision{These properties are instrumental for establishing
the strongly polynomial resolution of the QP (\ref{eq:QP}) for a matrix $M$ in this class; see Section~\ref{sec:irreducible Sbar}}.

\begin{proposition} \label{pr:irreducible barpsd} \rm
Let $A \in \overline{\mathbb{S}}_+^{\, n}$ be irreducible.  The following statements hold:

\gap

(a) If $d > 0$ satisfies $\overline{A} d \geq 0$, then either $\overline{A} d = 0$ or $\overline{A}$ (and thus $A$) is
positive definite.

\gap

(b) Every proper principal submatrix of $\overline{A}$, and thus of $A$, is positive definite.

\gap

(c) $( A/A_{\alpha \alpha} )_{ii} > 0$ for all $\alpha \subset [ n ]$ with
$| \alpha | \leq n - 2$ and all $i \not\in \alpha$.

\gap

(d) The kernel of $\overline{A}$ is of dimension at most one; moreover, if its dimension is
one, then there exists $f > 0$ such that all vectors in $\mbox{ker}(A)$ are multiples
of $f$.
\end{proposition}

\begin{proof} Since $\overline{A}$ is symmetric positive semidefinite, its irreducibility
implies that it has positive diagonals.  For (a), suppose that the vector $\overline{A} d$ contains
a positive element.  Then $\overline{A}$ is  irreducibly diagonally dominant; hence by \cite[Theorem~6.2.6]{Ortega90},
it follows that $\overline{A}$ is nonsingular, thus positive definite.  Hence so is $A$.

\gap

For (b), we may assume that $\overline{A}$ is singular.  Let $d > 0$ be such that $\overline{A}d = 0$.
\revision{This is equivalent to $\mbox{Diag}(d) \overline{A} \, \mbox{Diag}(d) \onebld = 0$,
where $\mbox{Diag}(d)$ is the diagonal matrix whose diagonal elements are the respective components
of the vector $d$ and $\onebld$ is the vector of all ones.  The matrix $\mbox{Diag}(d) \overline{A} \, \mbox{Diag}(d)$ remains
irreducible and belongs to $\overline{\mathbb{S}}_+^{\, n}$.}  Hence,
without loss of generality, we may take $d$ to be the vector of all ones from the start.  Thus, we have
\[
a_{jj} \, = \, \displaystyle{
\sum_{k \neq j}
} \, | \, a_{jk} \, | \epc \forall \, j \, \in \, [ n ].
\]
The claim about the positive definiteness of the proper principal
submatrices amounts to showing that for all $i \in [ n ]$,
\begin{equation} \label{eq:proof of proper}
[ \, x^{\top} \overline{A}x = 0 \mbox{ and } x_i \, = \, 0 \, ] \ \Rightarrow \ x \, = \, 0.
\end{equation}
Let $x$ be a vector satisfying the left-hand conditions, we have
\[
0 \, = \, x^{\top} \overline{A} x \, = \, \displaystyle{
\sum_{j \neq k}
} \, | \, a_{kj} \, | \, ( \, x_j - x_k \, )^2 + \displaystyle{
\sum_{j=1}^n
} \, \left( \, a _{jj} - \displaystyle{
\sum_{k \neq j}
} \, | \, a_{kj} \, | \, \right) x_j^2 \, = \, \displaystyle{
\sum_{j \neq k}
} \, | \, a_{kj} \, | \, ( \, x_j - x_k \, )^2.
\]
Thus, $| \, a_{kj} \, | \, ( \, x_j - x_k \, ) = 0$ for all $j \neq k$.  Since $A$ is
irreducible, for every index $k \neq i$,
there exists a sequence of indices $\{ i = i_0, i_1, \cdots, i_m = k \}$ such that
$a_{i \, i+1} \neq 0$ for all $j = 0, \cdots, m-1$.  Hence it follows that
\[
0 \, = \, x_i \, = \, x_{i_1} \, = \, \cdots \, = \, x_{i_m} \, = \, x_k.
\]
Thus (\ref{eq:proof of proper}) holds, establishing that every proper principal
submatrix of $\overline{A}$ is positive definite; thus the same holds for every proper
principal submatrix of $M$.

\gap

For (c), with $\alpha$ and $i$ as given,
$( A/A_{\alpha \alpha} )_{ii}$ is the Schur complement
of the positive definite matrix $A_{\alpha \alpha}$ in the positive definite matrix
\[
\left[ \begin{array}{ll}
A_{\alpha \alpha} & A_{\alpha i} \\ [0.1in]
A_{i \alpha} & a_{ii}	
\end{array} \right],
\]
thus the positivity of $( A/A_{\alpha \alpha} )_{ii}$ follows readily from the
Schur determinantal formula.

\gap

Finally, to prove (d), let $x$ satisfy $\overline{A} x = 0$ and $d > 0$ satisfy $\overline{A} d = 0$.
Following the proof of part (b), we can deduce that $\displaystyle{
\frac{x_i}{d_i}
} = \, \displaystyle{
\frac{x_k}{d_k}
}$ for all pairs of indices $j$ and $k$.  This easily yields the desired conclusion of (d).  Indeed,
if $x_k \neq 0$ for an index $k$, then $x_i = \displaystyle{
\frac{d_i}{d_k}
}\, x_k$ for all $i$ so that the desired vector $f$ has its components equal to
$\displaystyle{
\frac{d_i}{d_k}
}$ for all $i \in [ n ]$.
\end{proof}

Based on the above proposition, if $A \in \overline{\mathbb{S}}_+^{\, n}$ is irreducible, then a vector $d > 0$ satisfying
$\overline{A} d \geq 0$ can be obtained very easily by solving the system of linear equations $\overline{A}d = 0$ with $d_n = 1$.
Indeed, if this system has no solution, then $\overline{A}$ must be nonsingular.  In this case,
$d \triangleq \overline{A}^{\, -1} r$, where $r$ is any positive vector, is a desired vector.  If the system has a solution, then
the solution must be positive by part (d) of the proposition.

\section{The PPPA with $p$ satisfying (\ref{eq:psd n-step vector})} \label{sec:PPPA with p}

In the section, we assume that $M \in \mathbb{S}_+^n$ and that a vector $p$ exists
satisfying (\ref{eq:psd n-step vector}) and $q + \tau_0 p \geq 0$ for some scalar $\tau_0 > 0$.
Not required to be nonnegative, the vector $p$ satisfying the latter condition is sufficient to
initiate the PPPA algorithm for solving the equivalent linear complementarity formulation
of the QP (\ref{eq:QP}).  The goal of this section is to present a streamlined implementation of
the  algorithm that yields the $2n$-step Algorithm~I when the stronger assumption
(\ref{eq:n-step vector}) is in place; the discussion will elucidate the role played by the
condition (\ref{eq:psd n-step vector}) and prepare the proof for the same linear-step
termination when $M \in \overline{\mathbb{S}}_+^{\, n}$ is irreducible.  Generalizations
beyond the latter class will the the subject of Section~\ref{sec:extended problem classes}.

\subsection{Analysis of the pivot steps} \label{subsec:pivots}

\revision{The goal of the analysis here is to demonstrate two main things that are the key to the complexity
results in the next section:}

\gap

\revision{(A) after an $x$-variable reaches its upper bound, it stays at the upper bound
till the end of the algorithm;}

\gap

\revision{(B) a $2 \times 2$ pivot is always triggered by a nonbasic $w$-variable with a zero Schur diagonal entry,
where $w$ is the complement of $x$.}

\gap

\revision{We will return to these two points at the end of the analysis, which we will incorporate in the statement of
streamlined Algorithm~II.}

\gap

Suppose that a solution $x^{\rm cur}$ is available to (\ref{eq:parametric QP})
corresponding to a value $\tau_{\rm cur} > 0$.  Define three index sets
\begin{equation} \label{eq:3 index sets}
\begin{array}{ll}
\alpha_{\rm cur} \, \triangleq \, \left\{ \, i \, \in \, [ n ] \, \mid \,
0 \, < \, x^{\rm cur}_i \, < u_i \, \right\} & \mbox{basic nondegenerate; between bounds}
\\ [0.1in]
\beta_{\rm cur} \, \triangleq \, \left\{ \, i \, \in \, [ n ] \, \mid \,
x^{\rm cur}_i \, = \, 0 \, \right\} & \mbox{nonbasic at lower bounds} \\ [0.1in]
\gamma_{\rm cur} \, \triangleq \, \left\{ \, i \, \in \, [ n ] \, \mid \, u_i < \infty \mbox{ and }
x^{\rm cur}_i \, = \, u_i \, \right\} & \mbox{nonbasic at upper bounds}.
\end{array} \end{equation}
We wish to compute a solution path of (\ref{eq:parametric QP}) by decreasing $\tau$
until $\tau = 0$.
The PPPA accomplishes this computation by a sequence of simple diagonal pivots or
$2 \times 2$ pivots
if the former is not possible.  With (\ref{eq:3 index sets}) in the backdrop,
let $( \alpha,\beta,\gamma)$ be three index sets partitioning $[ n ]$.  At the current
value of $\tau$, the variables in these index sets are, respectively: between
bounds ($\alpha$), at lower bounds ($\beta$), and at upper bounds ($\gamma$).
We may write the
Karush-Kuhn-Tucker (KKT) conditions of the parametric QP (\ref{eq:parametric QP})
as follows, where it is understood that if an upper bound $u_i$ is infinite, the
corresponding slack $s_i$ and multiplier $\lambda_i$ are non-existent:
\begin{equation} \label{eq:ub-KKT}
\begin{array}{lll}
0 \, \leq \, w_{\alpha} \, = \, q_{\alpha} + \tau \, p_{\alpha} +
M_{\alpha \alpha} x_{\alpha} + M_{\alpha \beta} x_{\beta} + M_{\alpha \gamma} x_{\gamma}
+ \lambda_{\alpha} & \perp & x_{\alpha} \, \geq 0 \\ [0.1in]
0 \, \leq \, w_{\beta} \, = \, q_{\beta} + \tau \, p_{\beta} +
M_{\beta \alpha} x_{\alpha} + M_{\beta \beta} x_{\beta} + M_{\beta \gamma} x_{\gamma}
\hspace{0.3in} + \lambda_{\beta} & \perp & x_{\beta} \, \geq 0 \\ [0.1in]
0 \, \leq \, w_{\gamma} \, = \, q_{\gamma} + \tau \, p_{\gamma} +
M_{\gamma \alpha} x_{\alpha} + M_{\gamma \beta} x_{\beta} + M_{\gamma \gamma} x_{\gamma}
\hspace{0.6in} + \lambda_{\gamma} & \perp & x_{\gamma} \, \geq 0 \\ [0.1in]
0 \, \leq \, s_{\alpha} \, = \, u_{\alpha} \hspace{0.5in} - \hspace{0.2in}
x_{\alpha} & \perp & \lambda_{\alpha} \, \geq \, 0 \\ [0.1in]
0 \, \leq \, s_{\beta} \, = \, u_{\beta} \hspace{1.2in} - \hspace{0.1in}
x_{\beta} & \perp & \lambda_{\beta} \, \geq \, 0 \\ [0.1in]
0 \, \leq \, s_{\gamma} \, = \, u_{\gamma} \hspace{1.9in} - \hspace{0.1in}
x_{\gamma} & \perp & \lambda_{\gamma} \, \geq \, 0,
\end{array}
\end{equation}
where $\perp$ is the perpendicularity notation which in this context denotes the complementarity condition.
Assuming that $M_{\alpha \alpha}$ is nonsingular, we may perform a block pivot to make
the two tuples
\[
( x_{\alpha},w_{\beta},\lambda_{\gamma},s_{\alpha},s_{\beta},x_{\gamma} )
\epc \mbox{and} \epc
( w_{\alpha},x_{\beta},s_{\gamma},\lambda_{\alpha},\lambda_{\beta},w_{\gamma} )
\]
basic and nonbasic, respectively; this results in the equivalent system:
\begin{equation} \label{eq:basic-nonbasic KKT}
\begin{array}{l}
x_{\alpha} \, = \, -\bar{q}_{\alpha} - \tau \, \bar{p}_{\alpha} +
M_{\alpha \alpha}^{-1} \, w_{\alpha} -
M_{\alpha \alpha}^{-1} M_{\alpha \beta} x_{\beta}
+ M_{\alpha \alpha}^{-1} M_{\alpha \gamma} s_{\gamma}
- M_{\alpha \alpha}^{-1} \lambda_{\alpha} \\ [0.1in]
0 \, \leq \, x_{\alpha} \, \perp \, w_{\alpha} \, \geq \, 0 \\ [0.1in]
w_{\beta} \, = \, \bar{q}_{\beta} + \tau \, \bar{p}_{\beta} +
M_{\beta \alpha} \, M_{\alpha \alpha}^{-1} w_{\alpha} +
( M/M_{\alpha \alpha} )_{\beta \beta} \, x_{\beta} -
( M/M_{\alpha \alpha} )_{\beta \gamma} \, s_{\gamma}
- M_{\beta \alpha} M_{\alpha \alpha}^{-1} \lambda_{\alpha} +
\lambda_{\beta} \\ [0.1in]
0 \, \leq \, w_{\beta} \, \perp \, x_{\beta} \, \geq \, 0 \\ [0.1in]
\lambda_{\gamma} \, = \, -\bar{q}_{\gamma} - \tau \, \bar{p}_{\gamma}
- M_{\gamma \alpha} \, M_{\alpha \alpha}^{-1} w_{\alpha}
- ( M/M_{\alpha \alpha} )_{\gamma \beta} \, x_{\beta}
+ ( M/M_{\alpha \alpha} )_{\gamma \gamma} \, s_{\gamma}
+ M_{\gamma \alpha} M_{\alpha \alpha}^{-1} \lambda_{\alpha} + w_{\gamma} \\ [0.1in]
0 \, \leq \, \lambda_{\gamma} \, \perp \, s_{\gamma} \, \geq \, 0 \\ [0.1in]
s_{\alpha} \, = \, ( \, u_{\alpha} + \bar{q}_{\alpha} \, )
+ \tau \, \bar{p}_{\alpha} - M_{\alpha \alpha}^{-1} w_{\alpha} +
M_{\alpha \alpha}^{-1} M_{\alpha \beta} x_{\beta}
- M_{\alpha \alpha}^{-1} M_{\alpha \gamma} s_{\gamma}
+ M_{\alpha \alpha}^{-1} \, \lambda_{\alpha} \\ [0.1in]
0 \, \leq \, s_{\alpha} \, \perp \lambda_{\alpha} \, \geq \, 0 \\ [0.1in]
s_{\beta} \, = \, u_{\beta} - x_{\beta}; \epc
0 \, \leq \, s_{\beta} \, \perp \, \lambda_{\beta} \, \geq 0 \\ [0.1in]
x_{\gamma} \, = \, u_{\gamma} - s_{\gamma}; \epc
0 \, \leq \, x_{\gamma} \, \perp \, w_{\gamma} \, \geq 0,
\end{array} \end{equation}
where the updated vectors $\bar{q}$ and $\bar{p}$ are obtained from (\ref{eq:barpq-basic}) and (\ref{eq:barpq-nonbasic}).
By condition (\ref{eq:psd n-step vector}),
the components $\bar{p}_{\alpha} \geq 0$ and
$\bar{p}_{\gamma} \geq 0$; therefore, the basic variables $x_{\alpha}$ and
$\lambda_{\gamma}$ will increase in value when $\tau$ decreases before the next pivot.
Since the basic variables are kept nonnegative while the nonbasic variables are
set equal to zero throughout the PPPA, it follows that
$-\bar{q}_{\alpha} \geq 0$ and $-\bar{q}_{\gamma} \geq 0$.  At the current $\tau_{\rm cur}$, we have
\begin{equation} \label{eq:current nonnegative}
\left( \ -\bar{q}_{\alpha} - \tau_{\rm cur} \ \bar{p}_{\alpha}, \ \bar{q}_{\beta} + \tau_{\rm cur} \, \bar{p}_{\beta}, \
-\bar{q}_{\gamma} - \tau_{\rm cur} \, \bar{p}_{\gamma}, \  u_{\alpha} + \bar{q}_{\alpha} + \tau_{\rm cur} \, \bar{p}_{\alpha} \ \right) \, \geq 0.
\end{equation}
To determine the next critical value of
$\tau$ and thus the next pivot, we perform the ratio test:
\begin{equation} \label{eq:new critical value}
\tau_{\rm new} \, \triangleq \, \max\left\{ \, \underbrace{\displaystyle{
\max_{i \in \beta}
} \, \left\{ \, -\displaystyle{
\frac{\bar{q}_i}{\bar{p}_i}
} \, \mid \, \bar{p}_i > 0 \, \right\}}_{\mbox{denoted $\rho_0^{\, \tau}$}}; \
\underbrace{\displaystyle{
\max_{i \in \alpha}
} \, \left\{ \, -\displaystyle{
\frac{u_i + \bar{q}_i}{\bar{p}_i}
} \, \mid \, \bar{p}_i > 0 \mbox{ and } u_i < \infty \, \right\}}_{\mbox{denoted
$\rho_u^{\, \tau}$}}, \, 0, \, \right\}
\end{equation}
If $\tau_{\rm new} = 0$ (this is the same as (\ref{eq:ratio test in pd})),
then the current system (\ref{eq:basic-nonbasic KKT}) with
$\tau = 0$ yields an optimal solution to the original QP (\ref{eq:QP}).  So in what
follows, we assume without loss of generality that $\tau_{\rm new} > 0$.
There are two main cases to consider:

\gap

{\bf (1)} $\tau_{\rm new} = \rho_u^{\, \tau} = -\displaystyle{
\frac{u_{\bar{i}} + \bar{q}_{\bar{i}}}{\bar{p}_{\bar{i}}}
}$ for a maximizing index $\bar{i} \in \alpha$.

\gap

{\bf (2)} $\tau_{\rm new} = \rho_0^{\, \tau} = -\displaystyle{
\frac{\bar{q}_{\bar{i}}}{\bar{p}_{\bar{i}}}
}$ for a maximizing index $\bar{i} \in \beta$.  There are two subcases:
(2a) $( M/M_{\alpha \alpha} )_{\bar{i} \, \bar{i}} > 0$, or
(2b) $( M/M_{\alpha \alpha} )_{\bar{i} \, \bar{i}} = 0$.

\gap

Cases (1) and (2a) lead to the following updates of the index sets, respectively:
\begin{equation} \label{eq:index set 1}
\alpha_{\rm new} \, = \, \alpha \, \setminus \, \left\{ \, \bar{i} \, \right\}
\epc \mbox{and} \epc
\gamma_{\rm new} \, = \, \gamma \, \cup \, \left\{ \, \bar{i} \, \right\};
\end{equation}
%
%
%
\begin{equation} \label{eq:index set 2a}
\alpha_{\rm new} \, = \, \alpha \, \cup \, \left\{ \, \bar{i} \, \right\}
\epc \mbox{and} \epc
\beta_{\rm new} \, = \, \beta \, \setminus \, \left\{ \, \bar{i} \, \right\}.
\end{equation}
These updates are the same as those in Algorithm~I.

\gap

Case (2b) is the departure point from Algorithm~I and necessitates
a $2 \times 2$ pivot that exploits the symmetry and positive semidefiniteness of the Schur complement
$( M/M_{\alpha \alpha} )$.  In fact, $( M/M_{\alpha \alpha} )_{\bar{i} \, \bar{i}} = 0$ implies
that the $\bar{i}$th row and column of $( M/M_{\alpha \alpha} )$ are both entirely zero.
A second ratio test is needed to determine the outgoing basic variable, if there is one.
Letting
\[
\wh{M}_{\alpha \bar{i}} \, \triangleq \, M_{\alpha \alpha}^{-1}M_{\alpha \bar{i}}
\, = \,
\left[ \, M_{\bar{i} \alpha} M_{\alpha \alpha}^{-1} \, \right]^{\top}
\]
we perform the ratio test:
\begin{equation} \label{eq:another ratio}
\rho_{\min} \, \triangleq \, \min\left\{ \, \begin{array}{l}
\rho_>^{\, x} \, \triangleq \, \min\left\{ \displaystyle{
\frac{-\bar{q}_j - \tau_{\rm new} \, \bar{p}_j}{\wh{m}_{j \bar{i}}}
} \, \mid \, j \, \in \, \alpha, \ \wh{m}_{j \bar{i}} \, > \, 0 \, \right\}; \\ [0.25in]
\rho_<^{\, x} \, \triangleq \, \min\left\{\displaystyle{
\frac{u_j + \bar{q}_j + \tau_{\rm new} \, \bar{p}_j}{-\wh{m}_{j \bar{i}}}
} \, \mid \, j \, \in \, \alpha, \ u_j < \infty \ \mbox{ and } \
\wh{m}_{j \bar{i}} \, < \, 0 \, \right\}; \\ [0.2in]
\ u_{\bar{i}} \epc \mbox{if finite}
\end{array} \right\}.
\end{equation}
If there is no candidate qualified in the above ratio test, in particular if
$u_{\bar{i}} = \infty$, then the right-hand side of the $w_{\bar{i}}$-row in
the system (\ref{eq:basic-nonbasic KKT}) is negative for $\tau < \tau_{\rm new}$,
in particular for $\tau = 0$, contradicting the nonnegativity of $w_{\bar{i}}$.
Thus the QP (\ref{eq:QP}) has no solution when $\rho_{\min} = \infty$.  With a finite $\rho_{\min}$,
there are 3 cases to consider, depending which quantity yields the minimum
ratio $\rho_{\min}$.

\gap

$\bullet $ $\rho_{\min} = u_{\bar{i}}$: in this case, the change of the membership
of the index $\bar{i}$ is as follows:
\begin{equation} \label{eq:index set 2b_u}
\beta_{\rm new} \, = \, \beta \, \setminus \, \left\{ \, \bar{i} \, \right\}
\epc \mbox{and} \epc
\gamma_{\rm new} \, = \, \gamma \, \cup \, \left\{ \, \bar{i} \, \right\},
\end{equation}
implying an increase of one element in the union
$\alpha_{\rm new} \cup \gamma_{\rm new}$.  This update of the index sets corresponds
to a $2 \times 2$ pivot whereby the role of the presently basic pair
of variables $( w_{\bar{i}},s_{\bar{i}} )$ is exchanged with the presently nonbasic pair
$( x_{\bar{i}},\lambda_{\bar{i}} )$, making the former pair nonbasic and the latter pair
basic.  There is no matrix operations needed.  A new iteration commences with
the updated index sets (\ref{eq:index set 2b_u}).

\gap

$\bullet $ $\rho_{\min} = \rho^{\, x}_> = \displaystyle{
\frac{-\bar{q}_{\bar{j}} - \tau_{\rm new} \, \bar{p}_{\bar{j}}}{
\bar{m}_{\bar{j} \, \bar{i}}}
}$ for a minimizing index $\bar{j} \in \alpha$ with $\bar{m}_{\bar{j} \, \bar{i}} > 0$:
in this case, the change of the membership of the pair of indices
$( \, \bar{i}, \, \bar{j} \, )$ is as follows:
\begin{equation} \label{eq:index set 2b_>}
\beta_{\rm new} \, = \, \beta \, \setminus \, \left\{ \, \bar{i} \, \right\} \, \cup
\, \left\{ \, \bar{j} \, \right\}
\epc \mbox{and} \epc
\alpha_{\rm new} \, = \, \alpha \, \cup \, \left\{ \, \bar{i} \, \right\} \, \setminus
\, \left\{ \, \bar{j} \, \right\},
\end{equation}
implying no change in cardinality of the set $\alpha_{\rm new}$ and no change in
$\gamma$.  This update of the index sets corresponds to the $2 \times 2$ pivot on the
submatrix: $\left[ \begin{array}{cc}
\bar{m}_{\bar{j} \, \bar{j}} & - \bar{m}_{\bar{j} \, \bar{i}} \\
\bar{m}_{\bar{i} \, \bar{j}} & 0	
\end{array} \right]$, where $\bar{m}_{\bar{j} \, \bar{j}}$ is the $\bar{j}$-diagonal
entry of $M_{\alpha \alpha}^{-1}$.

\gap

$\bullet $ $\rho_{\min} = \rho^{\, x}_< = \displaystyle{
\frac{u_{\bar{j}} + \bar{q}_{\bar{j}} + \tau_{\rm new} \, \bar{p}_{\bar{j}}}{
-\bar{m}_{\bar{j} \, \bar{i}}}
}$ for a minimizing index $\bar{j} \in \alpha$ with $u_{\bar{j}} < \infty$ and
$\bar{m}_{\bar{j} \, \bar{i}} < 0$:
this is the mirror of the previous case with the following change of the membership of the pair
of indices $( \, \bar{j}, \, \bar{i} \, )$:
\begin{equation} \label{eq:index set 2b_<}
\beta_{\rm new} \, = \, \beta \, \setminus \, \left\{ \, \bar{i} \, \right\}; \epc
\alpha_{\rm new} \, = \, \alpha \, \cup \, \left\{ \, \bar{i} \, \right\} \, \setminus
\, \left\{ \, \bar{j} \, \right\};
\epc \mbox{and} \epc
\gamma_{\rm new} \, = \, \gamma \, \cup \, \left\{ \, \bar{j} \, \right\}.
\end{equation}
The difference from (\ref{eq:index set 2b_>}) is that the index $\bar{j} \in \alpha$
is transitioned to $\gamma$ instead of to $\beta$.
This update of the index sets corresponds to the $2 \times 2$ pivot on the
submatrix: $\left[ \begin{array}{cc}
\bar{m}_{\bar{j} \, \bar{j}} & \bar{m}_{\bar{j} \, \bar{i}} \\
-\bar{m}_{\bar{i} \, \bar{j}} & 0	
\end{array} \right]$.

\gap

\revision{The above analysis supports the two points (A) and (B) made at the beginning of this subsection:}

\gap

\revision{(A) can be seen from the updates of the pivot indices:
(\ref{eq:index set 1}), (\ref{eq:index set 2a}), (\ref{eq:index set 2b_u}), and
(\ref{eq:index set 2b_>}), and (\ref{eq:index set 2b_<}) that show no index leaves
the index set $\gamma$ once it is there;}

\gap

\revision{(B) happens when $( M/M_{\alpha \alpha} )_{\bar{i} \bar{i}} = 0$ where $w_{\bar{i}}$ with
$\bar{i} \in \beta$ is the nonbasic variable to be pivoted on.  The resulting $2 \times 2$ pivot
involves the change of role of $\bar{i}$ with an index $\bar{j} \in \alpha$;
$\bar{i}$ will always move to the (new) index set $\alpha$, where $\bar{j}$ will move either
to the (new) $\beta$ as in (\ref{eq:index set 2b_>}) or to the (new) $\gamma$ as in (\ref{eq:index set 2b_<}).}

\gap

While not used in the subsequent analysis, it is interesting to note that the index set $\beta$ will not grow.


\gap

\rule{6in}{0.01in}

{\bf Algorithm II: The streamlined PPPA for (\ref{eq:QP}) with $M$ positive semidefinite}

\rule{6in}{0.01in}

\gap

\noindent {\bf Input.} Let $M \in \mathbb{S}_+^n$ have positive diagonal
entries.  Let $p \in \mathbb{R}^n$
satisfy (\ref{eq:psd n-step vector}) and such that $q + \tau_0 p \geq 0$ for some $\tau_0 > 0$.

\gap

\noindent {\bf Initialization.}  Let $\alpha = \gamma = \emptyset$ and $\beta = [ n ]$.

\gap

\noindent {\bf Main computations.}  Compute the vectors $\bar{q}$ and $\bar{p}$ by
(\ref{eq:barpq-basic}) and (\ref{eq:barpq-nonbasic}) and determine the new critical value $\tau_{\rm new}$
by performing the ratio test (\ref{eq:new critical value}).   If $\tau_{\rm new} = 0$; an optimal
solution of (\ref{eq:QP}) is on hand and given by $( x_{\alpha},x_{\beta},x_{\gamma} ) =
( -\bar{q}_{\alpha}, 0, u_{\gamma} )$.  Otherwise, proceed to the next step.
\gap

\noindent {\bf Updates of index sets.} Update the index triplet
$( \alpha,\beta,\gamma )$ according to (\ref{eq:index set 1}), (\ref{eq:index set 2a}),
(\ref{eq:index set 2b_u}),
(\ref{eq:index set 2b_>}), or (\ref{eq:index set 2b_<}) depending on the pivot
determinations.  In particular, if $\rho_{\min} = \infty$, then the QP (\ref{eq:QP})
is unbounded below.
Otherwise, return to the main computational step to continue the pivots.
\hfill $\Box$

\gap

\rule{6in}{0.01in}

\section{The Case of a $\overline{\mathbb{S}}^{\, n}_+$-Matrix} \label{sec:irreducible Sbar}

Let $M$ in $\overline{\mathbb{S}}_+^{\, n}$ be irreducible.  Without loss of generality, we may assume
$d \in \mathbb{R}^n_{++}$ exists such that $\overline{M}d = 0$.  Let $p \triangleq Md$.  Then the vector $p$ is nonnegative
and satisfies (\ref{eq:psd n-step vector}).  Nevertheless, there may not exist $\tau_0$ such that $q + \tau_0 p \geq 0$,
if for instance $p$ has a zero component and the corresponding $q$ component is negative.  We postpone the treatment of
this detail until later.  The theorem below assumes that the PPPA can be initiated with the vector $p$ as given.

\begin{theorem} \label{th:irreducible barpsd} \rm
Let $M \in \bar{\mathbb{S}}_+^{\, n}$ be irreducible.  Let $d > 0$ satisfy $\overline{M}d = 0$.
Suppose $q + \tau_0 p \geq 0$ for some $\tau_0 > 0$, where $p \triangleq Md$.  Then in at most $2n$ pivots
of the kinds (\ref{eq:index set 1}) or (\ref{eq:index set 2a}), the PPPA will either compute an optimal
solution of the QP (\ref{eq:QP}) or certify that the QP has no finite optimal solution.  Hence,
the overall complexity of the PPPA is of order O$(n^3)$.
\end{theorem}

\begin{proof}  By Proposition~\ref{pr:irreducible barpsd}(c), it follows that all pivots of the PPPA are diagonal
until $| \alpha | = n - 1$.
Moreover throughout these pivots, a between-bound $x$-variable can only reach its upper bound; when that happens, it stays at the upper bound.
Now consider the case where in the next pivot, the only nonbasic $x$-variable, say $x_{\bar{i}}$ for some $\bar{i} \in \beta$, is the candidate
to be pivoted on and
the Schur complement $m_{\bar{i} \bar{i}} - M_{\bar{i} \alpha} M_{\alpha \alpha}^{-1} M_{\alpha \bar{i}} = 0$.  By part (b) of the same proposition,
the $\gamma$ set at this point must be empty; otherwise the Schur entry would be positive.  This triggers a $2 \times 2$ pivot
of the kind (\ref{eq:index set 2b_u}), (\ref{eq:index set 2b_>}), or (\ref{eq:index set 2b_<}), unless $\rho_{\min} = \infty$ in which case
the algorithm terminates, confirming no finite optimal solution of the QP (\ref{eq:QP}).  The first and third kind will empty the set
$\beta$.  From that point on, the remaining pivots are of type (\ref{eq:index set 1}) each driving a between-bound $x$-variable to its upper bound.  If
the pivot is of the second kind with a finite $\rho_{\min}$, i.e., (\ref{eq:index set 2b_>}) occurs, then the index $\bar{i} \in \beta$ is exchanged
with an index $\bar{j} \in \alpha$.
It can be shown, by the properties of the $2 \times 2$ pivots,
that $\bar{j}$ will stay in the new $\beta$ in the immediately next pivot, if the algorithm does not yet terminate.  That next pivot
must be of the kind (\ref{eq:index set 1}), rendering the $\gamma$ set nonempty.  Since $| \gamma |$, which is now nonempty, can only increase,
the diagonal entry corresponding to the $x_{\bar{j}}$-variable will be positive in all subsequent pivots.  Consequently, in the remaining pivots
which are all of type (\ref{eq:index set 1}),
either $x_{\bar{j}}$ remains the only
nonbasic $x$-variable till the end, i.e., $\bar{j}$ is the only element of the $\beta$ set, or the $\beta$-set will be emptied out.  When the latter occurs,
the only pivots are the transitions of the between-bound $x$-variables reaching their respective upper bounds.  The claim of $2n$-step termination
of the theorem now follows readily, and so does the overall O$(n^3)$ complexity of the algorithm.
\end{proof}

We next address the technical issue where the vector $p$ in Theorem~\ref{th:irreducible barpsd}
is such that $p_i = 0$ and $q_i < 0$ for some index $i$ so that the PPPA cannot
be initiated.  Without assuming irreducibility, the
following lemma shows that in this case, the original QP (\ref{eq:QP}) can be reduced in
a way depending on whether $u_i$ is infinite or finite.  In both parts of the lemma, equivalence means that there is a one-to-one
correspondence between the optimal solutions of the original QP and the
respective reduced QPs.

\begin{lemma} \label{lm:nonpositive row} \rm
Let $M \in \overline{\mathbb{S}}_+^{\, n}$.  Let $d > 0$ satisfy $\overline{M}d \geq 0$
and let $p \triangleq \thalf ( M + \overline{M} )d$.  Suppose $p_i = 0$, $m_{ii} > 0$, and $q_i < 0$.
Then the following two
statements hold.

\gap

(a) If $u_i = \infty$, then (\ref{eq:QP}) is equivalent to
\begin{equation} \label{eq:infinite ub reduced QP}
\begin{array}{ll}
\displaystyle{
\operatornamewithlimits{\mbox{\bf minimize}}_{0 \, \leq \, x^{-i} \, \leq \, u^{-i}}
} & \wh{q}^{\, \top}x^{-i} + \thalf \, ( x^{-i} )^{\top} ( M/m_{ii} ) x^{-i},
\end{array} \end{equation}
where $u^{-i}$ is the upper-bound vector $u$ without the $i$th component and $\wh{q}$ is
the principal pivot transform of the vector $q$ without the $i$th component, i.e.,
\begin{equation} \label{eq:pivot transform of q}
\wh{q}_j \, \triangleq \, q_j - \displaystyle{
\frac{m_{ji}}{m_{ii}}
} \, q_i, \epc j \, \in \, [ n ] \, \setminus \{ i \}.
\end{equation}
(b) If $u_i < \infty$, then (\ref{eq:QP}) is equivalent to
\begin{equation} \label{eq:finite ub reduced QP}
\begin{array}{ll}
\displaystyle{
\operatornamewithlimits{\mbox{\bf minimize}}_z
} & \wt{q}^{\, \top}z + \thalf \, z^{\top} \wt{M} z \\ [0.1in]
\mbox{\bf subject to} & 0 \, \leq \, z_j \, \leq \, u_j, \epc \forall \, j \, \neq \, i
\\ [5pt]
\mbox{\bf and} & 0 \, \leq \, z_i,
\end{array} \end{equation}
where $\wt{q} \in \mathbb{R}^n$ has components given by
\[
\wt{q}_j \, \triangleq \, \left\{ \begin{array}{cl}
q_j + m_{ji}u_i & \mbox{if $j \neq i$}	\\ [5pt]
-\left( \, q_i + m_{ii} u_i \, \right) & \mbox{if $j = i$},
\end{array} \right.
\]
and $\wt{M} \triangleq \mathbb{I}_{-1}^{\, i} M \, \mathbb{I}_{-1}^{\, i}$ with
$\mathbb{I}_{-1}^{\, i}$ being the identity matrix with the $i$th diagonal element
changed to $-1$.
\end{lemma}

\begin{proof} For notational simplicity, we assume that all upper bounds $u_j$ for
$j \neq i$ are finite.  Introducing multipliers $v_j$ for the upper-bound constraints,
We can write the Karush-Kuhn-Tucker (KKT) conditions of the QP (\ref{eq:QP}) as follows:
with $\perp$ denoting the perpendicularity notation which in the context describes
the complementary slackness condition:

\gap

$\bullet $ if $u_i = \infty$,
\begin{equation} \label{eq:KKT infinite ub}
\left\{ \begin{array}{lll}
0 \, \leq \, q_j + \displaystyle{
\sum_{k \neq i}
} \, m_{jk} \, x_k + m_{ji} \, x_i + v_j & \perp & x_j \, \geq \, 0 \\ [0.2in]
0 \, \leq \, u_j - x_j & \perp & v_j \, \geq \, 0
\end{array} \right\}  \epc \forall \, j \, \neq \, i
\end{equation}
and $0 \, \leq \, q_i + \displaystyle{
\sum_{k \neq i}
} \, m_{ik} \, x_k + m_{ii} \, x_i \, \perp \, x_i \, \geq \, 0$;

\gap

$\bullet $ if $u_i < \infty$,
\begin{equation} \label{eq:KKT finite ub}
\hspace{0.4in} \left\{ \begin{array}{lll}
0 \, \leq \, q_j + \displaystyle{
\sum_{k \neq i}
} \, m_{jk} \, x_k + m_{ji} \, x_i + v_j & \perp & x_j \, \geq \, 0 \\ [0.2in]
0 \, \leq \, u_j - x_j & \perp & v_j \, \geq \, 0
\end{array} \right\}  \epc \forall \, j \, = \, 1, \cdots, n.
\end{equation}
Consider first the case of an infinite $u_i$.  By
part (c) of Lemma~\ref{lm:nonpositive row},
we know that $m_{ik} \leq 0$ for all $k \neq i$.  Hence, we must have
\[
0 \, = \, q_i + \displaystyle{
\sum_{k \neq i}
} \, m_{ik} \, x_k + m_{ii} \, x_i,
\]
in any optimal solution of (\ref{eq:QP}), because if the right-hand expression were positive for
some such solution, then by complementarity, it must follow that
$x_i = 0$, which is a contradiction.  The above equation yields
\begin{equation} \label{eq:yielding xi}
x_i \, = \, \displaystyle{
\frac{-1}{m_{ii}}
} \, \left[ \, q_i + \displaystyle{
\sum_{k \neq i}
} \, m_{ik} \, x_k \, \right].
\end{equation}
Substituting this expression for $x_i$ into (\ref{eq:KKT infinite ub}) yields the KKT
system for the QP (\ref{eq:infinite ub reduced QP}).  This shows that every optimal solution
of (\ref{eq:QP}) yields an optimal solution of (\ref{eq:infinite ub reduced QP}) by
simply removing the $i$th variable.  Conversely,
for any optimal solution $x^{-i}$ of the latter QP without the $i$th variable, the
expression (\ref{eq:yielding xi}) defines a nonnegative $x^{\, i}$ which together with
$x^{-i}$ satisfies the KKT conditions of the QP (\ref{eq:QP}).  This completes the proof
of the equivalence of the two QPs (\ref{eq:QP}) and (\ref{eq:infinite ub reduced QP}).

\gap

Suppose $u_i < \infty$.  Then similarly, we have
\begin{equation} \label{eq:ith complementarity finite ub}
0 \, = \, q_i + \displaystyle{
\sum_{k \neq i}
} \, m_{ik} \, x_k + m_{ii} \, x_i + v_i,
\end{equation}
in any solution of the KKT system (\ref{eq:KKT finite ub}).
Letting $z_i \triangleq u_i - x_i$, or equivalent, $x_i = u_i - z_i$, and substituting
the latter expression for $x_i$ into the KKT conditions (\ref{eq:KKT finite ub}) for
all $j \neq i$, we obtain
\[
\left\{ \begin{array}{lll}
0 \, \leq \, \wt{q}_j + \displaystyle{
\sum_{k \neq i}
} \, m_{jk} \, x_k - m_{ji} \, z_i + v_j & \perp & x_j \, \geq \, 0 \\ [0.2in]
0 \, \leq \, u_j - x_j & \perp & v_j \, \geq \, 0
\end{array} \right\}  \epc \forall \, j \, \in \, [ n ] \setminus \{ i \}
\]
and
\[ \hspace{-1.5in}
0 \, \leq \, v_i \, = \, \wt{q}_i - \displaystyle{
\sum_{k \neq i}
} \, m_{ik} \, x_k + m_{ii} \, z_i \, \perp \, z_i \, \geq \, 0.
\]
This shows that every optimal solution
of (\ref{eq:QP}) yields an optimal solution of (\ref{eq:finite ub reduced QP}) with
all variables $x_j$ for $j \neq i$ remaining the same and the $i$th variable is the slack
of the $i$th upper-bound constraint.  Conversely, if $( z_i,x^{-i} )$ is an optimal
solution of (\ref{eq:finite ub reduced QP}), then with
$x_i \triangleq u_i - z_i$, (\ref{eq:ith complementarity finite ub}) holds and this
implies that $x_i \geq 0$.  To deduce the latter nonnegativity, suppose $x_i < 0$.
Then $v_i > 0$, which implies by complementarity, $z_i = 0$; but then $x_i = u_i > 0$,
which is a contradiction.  Hence $( x_i,x^{-i} )$ is an optimal solution of
(\ref{eq:QP}).
\end{proof}

By part (b) of Lemma~\ref{lm:schur psd}, the Schur complement $(M/m_{ii})$ is in
class $\overline{\mathbb{S}}_+^{\, n-1}$.  Thus in situation (a) of
Lemma~\ref{lm:nonpositive row}, we have successfully reduced the
original QP (\ref{eq:QP}) to a problem of the same type by one dimension less.  We can next
derive a desired parametric vector to the reduced problem
(\ref{eq:infinite ub reduced QP}).  In situation (b) of
Lemma~\ref{lm:nonpositive row}, while the matrix
$\wt{M}$ remains in class $\overline{\mathbb{S}}_+^{\, n}$, its dimension is not reduced.
Nevertheless, one variable is no longer upper bounded.  Thus when we next repeat the
solution strategy on the problem (\ref{eq:finite ub reduced QP}), we obtain an
equivalent problem either with one less variable or with one more variable not upper
bounded.  Continuing this way, if we fail
to obtain an equivalent problem of lower dimension, then we
will arrive, after $n$ occurrrences of the situation (b) of
Lemma~\ref{lm:nonpositive row}, at a problem of the same type as (\ref{eq:QP})
with all the
variables not upper bounded.  For such a problem without upper bounds, part (a) of the
Lemma is applicable to reduce the dimension of the problem if a desired vector $p$ still
has not been obtained. In summary, via the repeated application of
Lemma~\ref{lm:nonpositive row}, we will eventually obtain a problem, possibly of reduced
size, with a desired vector $p$ satisfying the initialization requirement of the PPPA.
Therefore, the QP (\ref{eq:QP}) with $M \in \overline{S}_+^{\, n}$ can be resolved in
strongly polynomial time by repeated reduction via Lemma~\ref{lm:nonpositive row}
until Theorem~\ref{th:irreducible barpsd} is applicable.

\subsection{Reducibility and tridiagonality} \label{subsec:reducibility and tridiagonality}

We complete the discussion of this section with two worthy points.  One is when
$M \in \overline{\mathbb{S}}_+^{\, n}$ is reducible and the other is when $M \in \overline{\mathbb{S}}_+^{\, n}$ is tridiagonal.
In the former case, $M$ can be
decomposed, through a principal rearrangement of its rows and same columns, into a block diagonal matrix with $L$ diagonal blocks
each being an irreducible matrix in $\overline{\mathbb{S}}_+^{\, m_{\ell}}$ with $\displaystyle{
\sum_{\ell=1}^L
} \, m_{\ell} = n$, for some positive integers $\{ \, m_{\ell} \, \}_{\ell=1}^L$.  Thus
by applying the PPPA to each of the decomposed quadratic subprograms, and invoking Lemma~\ref{lm:nonpositive row} if needed,
the overall complexity of resolving the QP (\ref{eq:QP}) is $\displaystyle{
\sum_{\ell=1}^L
} \, \mbox{O}(m_{\ell}^3)$, which is no bigger than $\mbox{O}(n^3)$.

\gap

When $M \in \overline{\mathbb{S}}_+^{\, n}$ is tridiagonal, the Cholesky factorization of the principal submatrices of $M$
can be obtained in linear time by successive pivoting on the diagonal entries of the Schur complements
with each pivot requiring only a constant number of operations (cf.\ the proof of part (a) in Proposition~\ref{pr:some further subclasses}).
Consequently, the algebraic computations in (\ref{eq:barpq-basic})
and (\ref{eq:barpq-nonbasic}) can be accomplished in linear time by such factorization.
As a result, the total complexity of solving the problem (\ref{eq:QP})
with a tridiagonal $M \in \mathbb{S}_+^n$ is O$(n^2)$.  At this time, it is not clear whether
there is a linear-time algorithm for this special class of problems.

\section{Some Extended Problem Classes} \label{sec:extended problem classes}

Building on Section~\ref{sec:irreducible Sbar} we consider extended problems where the matrix $M$ is ``$k$-weakly quasi-diagonally dominant''
for some (presumably small) integer $k \geq 1$.  One way to describe such a matrix is that it belongs to the class $\overline{\mathbb{S}}_+^{\, n;k}$;
that is, $M \in \mathbb{S}^n_+$ and all of its principal submatrices of order $n - k$ belong to $\overline{\mathbb{S}}_+^{\, n-k}$.
An equivalent definition is that for each \revision{$(n - k)$}-subset $\kappa$ of $[ n ]$, i.e., $\kappa$
contains \revision{$n - k$} indices, there exists a positive
vector $d^{\, \kappa} \in \mathbb{R}^{n - k}$ such that $\overline{M}_{\kappa \kappa} d^{\kappa} \geq 0$.
To see how the PPPA can be applied to QP (\ref{eq:QP}) with such a matrix $M$, we first consider the case $k = 1$.
With the aid of some post-PPPA lower-order matrix-vector operations, we obtain an algorithm with total complexity
$\mbox{O}( 2n \times n^3)$.  For a general $k$, this complexity becomes $\mbox{O}(( 2n )^k \times n^3)$,
which although remains strongly polynomial in $n$ for fixed $k$ practically may not be most desirable when $k$ exceeds a threshold.
In what follows, we explain the solution method for $k = 1$, briefly
sketch the case $k = 2$, and omit the details of the further generalizations as they all follow the same strategy for the first case.

\gap

Consider a matrix $M \in \overline{\mathbb{S}}_+^{\, n;1}$.  The development below follows the strategy described
in \cite[Section~5]{AdlerCottlePang16} for the P-case and extended herein to the positive semidefinite case.
  For each index $i \in [ n ]$, define two subproblems each of
order $n - 1$ obtained
by fixing the variable $x_i$ at either its lower bound or upper bound (if $u_i < \infty$):

\gap

$\bullet $ the lower-bound-reduced problem, labelled QP$_i^{\, \rm lb}$:
\begin{equation} \label{eq:fixed lb}
\begin{array}{ll}
\displaystyle{
\operatornamewithlimits{\mbox{\bf minimize}}_{0 \, \leq \, x^{\, -i} \, \leq \, u^{\, -i}}
} & ( q^{\, -i} )^{\, \top}x^{\, -i} + \thalf \, ( x^{\, -i} )^{\top} M^{\, -i} x^{\, -i},
\end{array} \end{equation}
where the superscript ``$-i$'' on $x$, $u$, and $q$ means the subvectors of $x$, $u$, and $q$ with their respective $i$th components removed,
and $M^{-i}$ is the principal submatrix of $M$ with the $i$th row and column removed;

\gap

$\bullet $ the upper-bound-reduced problem, labelled QP$_i^{\, \rm ub}$ (if $u_i < \infty$):
\begin{equation} \label{eq:fixed ub}
\begin{array}{ll}
\displaystyle{
\operatornamewithlimits{\mbox{\bf minimize}}_{0 \, \leq \, x^{\, -i} \, \leq \, u^{\, -i}}
} & ( q^{\, -i;\rm ub} )^{\, \top}x^{\, -i} + \thalf \, ( x^{\, -i} )^{\top} M^{\, -i} x^{\, -i},
\end{array} \end{equation}
where $q^{\, -i;\rm ub}_j \triangleq q_j + m_{ji} \, u_i$ for all $j \in [ n ] \setminus \{ \, i \, \}$.

\gap

The connections of these problems to the original QP (\ref{eq:QP}) are summarized in the result below.

\begin{proposition} \label{pr:reduction xi} \rm
The following three statements hold for any matrix $M \in \mathbb{S}^{\, n}$:

\gap

(i) If for some index $i \in [ n ]$, the QP$_i^{\, \rm lb}$ is unbounded below, then so is (\ref{eq:QP});
if QP$_i^{\, \rm lb}$ has an optimal solution $x^{\, -i;\rm lb}$ such that
\begin{equation} \label{eq:xi = 0}
q_i + \displaystyle{
\sum_{j \neq i}
} \, m_{ij} x^{\, -i;\rm lb}_j \, \geq \, 0,
\end{equation}
then the vector $x^{\, \rm lb} \triangleq ( 0, x^{\, -i;\rm lb} )$ is an optimal solution of (\ref{eq:QP}).

\gap

(ii) If for some index $i \in [ n ]$, the QP$_i^{\, \rm ub}$ is unbounded below, then so is (\ref{eq:QP});
if QP$_i^{\, \rm ub}$ has an optimal solution $x^{\, -i;\rm ub}$ such that
\begin{equation} \label{eq:xi = ui}
q_i + \displaystyle{
\sum_{j \neq i}
} \, m_{ij} x^{\, -i;\rm ub}_j \, \leq \, 0,
\end{equation}
then the vector $x^{\, \rm ub} \triangleq ( u_i, x^{\, -i;\rm ub} )$ is an optimal solution of (\ref{eq:QP}).

\gap

(iii) If for all $i \in [ n ]$, the QP$_i^{\, \rm lb}$ has no optimal solution $x^{\, -i;\rm lb}$ satisfying (\ref{eq:xi = 0}) and
the QP$_i^{\, \rm ub}$ has no optimal solution $x^{\, -i;\rm ub}$ satisfying (\ref{eq:xi = ui}), then either the QP (\ref{eq:QP})
has no optimal solution or every component of every one of its solution(s) is strictly between its bounds.
\end{proposition}

\begin{proof} Every feasible solution of either QP$_i^{\, \rm lb}$ or QP$_i^{\, \rm ub}$ readily yields a feasible solution of (\ref{eq:QP});
thus if either one of the two former QPs is unbounded below, then so is (\ref{eq:QP}).  The proof of
statements (i) and (ii) can be completed by verifying that the two vectors $x^{\, \rm lb}$ and $x^{\, \rm ub}$ satisfy the
KKT conditions of the QP (\ref{eq:QP}) if (\ref{eq:xi = 0}) and (\ref{eq:xi = ui}) hold, respectively.  Statement (iii) holds because
if $x^{\, \rm opt}$ is an optimal solution of
(\ref{eq:QP}), then for every $i \in [ n ]$, $x^{\, \rm opt}_i$ is either strictly between its two bounds or at one of them.
If $x^{\, \rm opt}_i$ is at its lower bound, then $x^{\, {\rm opt};-i}$ is an
optimal solution of QP$_i^{\, \rm lb}$ satisfying (\ref{eq:xi = 0}).  If $x^{\, \rm opt}_i$ is at its upper bound, then
$x^{\, {\rm opt};-i}$ is an optimal solution of QP$_i^{\, \rm ub}$ satisfying (\ref{eq:xi = ui}).
\end{proof}

We next address the question of how to verify for a given $i \in [ n ]$, if QP$_i^{\, \rm lb}$ (QP$_i^{\, \rm ub}$) has an optimal solution satisfying
(\ref{eq:xi = 0}) ((\ref{eq:xi = ui}), respectively).  We discuss only the former as the latter is the same.  To begin, note that with
$M \in \overline{\mathbb{S}}_+^{\, n;1}$, the PPPA can applied to QP$_i^{\, \rm lb}$.  Thus, in no more than $n - 1$ steps, either an optimal solution
to QP$_i^{\, \rm lb}$ is obtained, or it is determined that this problem has no optimal solution.  In the latter case, so does (\ref{eq:QP}).
In the former case, we claim that $\displaystyle{
\sum_{j \neq i}
} \, m_{ij} x^{\, -i;\rm lb}_j$ is a constant on the set of optimal solutions $x^{\, -i;\rm lb}$ of (\ref{eq:fixed lb}).  Indeed, if
$\wh{x}^{\, -i}$ and $\wt{x}^{\, -i}$ are any two solutions of (\ref{eq:fixed lb}), then $M^{\, -i} \, \wh{x}^{\, -i} = M^{\, -i}  \, \wt{x}^{\, -i}$,
by the symmetry and positive semidefinitness of $M^{\, -i}$.  This then implies $\displaystyle{
\sum_{j \neq i}
} \, m_{ij} \wh{x}^{\, -i}_j = \displaystyle{
\sum_{j \neq i}
} \, m_{ij} \wt{x}^{\, -i}_j$ by the symmetry and positive semidefinitness of $M$.  Consequently, with a single solution computed by the PPPA,
which has $\mbox{O}( ( n-1 )^3 )$-complexity, it can be determined
if (\ref{eq:fixed lb}) has an optimal solution $x^{\, -i;\rm lb}$ satisfying (\ref{eq:xi = 0}).

\gap

To check if the QP (\ref{eq:QP}) has an optimal solution in the case of statement (iii) in Proposition~\ref{pr:reduction xi},
it suffices to consider the system of linear inequalities:
\begin{equation} \label{eq:system of linear inequalities}
q + Mx \, = \, 0; \epc 0 \, \leq \, x \, \leq \, u.
\end{equation}
Indeed, if we are in the situation of statement (iii) of the proposition, then it is easy to see that the QP (\ref{eq:QP}) has an optimal solution
if and only if (\ref{eq:system of linear inequalities}) is consistent.  In turn, to check the latter, we may assume that $M$ is irreducible.  Then
with $M \in \overline{\mathbb{S}}_+^{\, n;1}$, let $\alpha \subset [ n ]$ be arbitrary with cardinality $n - 2$ and with complement $\beta$.  It then
follows that $x$ satisfies (\ref{eq:system of linear inequalities}) if and only if
\begin{equation} \label{eq:final-check 2var}
\begin{array}{lll}
0 & = & ( \, q_{\beta} - M_{\beta \alpha} M_{\alpha \alpha}^{-1} q_{\alpha} \, ) + ( M/M_{\alpha \alpha} ) x_{\beta}, \epc
0 \, \leq \, x_{\beta} \, \leq \, u_{\beta} \\ [0.1in]
0 & \leq & - M_{\alpha \alpha}^{-1} \left[ \, q_{\alpha} + M_{\alpha \beta} x_{\beta} \, \right] \, \leq \, u_{\alpha}.
\end{array} \end{equation}
The latter is a system of linear inequalities in 2 variables ($x_{\beta}$) with bounds, 2 equations, and $2 (n - 2)$ inequalities.  Thus its feasibility can
be resolved, say by the Fourier-Motzkin elimination method, in $\mbox{O}(n^2)$ operations.

\gap

Summarizing the above discussion, we present the following algorithm.

\gap

\rule{6in}{0.01in}

\gap

{\bf Algorithm~III: Resolving (\ref{eq:QP}) with $M \in \overline{\mathbb{S}}_+^{\, n;1}$}

\rule{6in}{0.01in}

\gap

\noindent {\bf Input.} Let $M \in \overline{\mathbb{S}}_+^{\, n;1}$ be given.

\gap

\noindent {\bf Main computations:}  For every $i \in [ n ]$, do the following:

\gap

$\bullet $ Apply PPPA to (\ref{eq:fixed lb}).  If this problem has no finite optimal solution,
then so does (\ref{eq:QP}); stop.  Otherwise, let ${x}^{\, -i;\rm lb}$ be a solution of (\ref{eq:fixed lb}) computed from PPPA.  If
(\ref{eq:xi = 0}) holds, then $x^{\, \rm lb} \triangleq ( 0,x^{\, -i;\rm lb} )$ is a desired solution of (\ref{eq:fixed lb}); stop.  Otherwise continue.

\gap

$\bullet $ If $u_i < \infty$, apply PPPA to (\ref{eq:fixed ub}). If this problem has no finite optimal solution, then so does (\ref{eq:QP}); stop.
Otherwise, let $x^{\, -i;\rm ub}$ be a solution of (\ref{eq:fixed ub}) computed from PPPA.   If
(\ref{eq:xi = ui}) holds, then $x^{\, \rm ub} \triangleq ( u_i,x^{\, -i;\rm ub} )$ is a desired solution of (\ref{eq:fixed lb}); stop.
Otherwise continue.

\gap

$\bullet $ Check the feasibility of the 2-variable system of linear inequalities (\ref{eq:final-check 2var}) by say the Fourier-Motzkin elimination method.
If it is not feasible, then the QP (\ref{eq:QP}) has no finite optimal solution; otherwise a solution is obtained.
\hfill $\Box$

\gap

\rule{6in}{0.01in}

\gap

\begin{theorem} \label{th:PPPA psd n-2 rank} \rm
Let $M \in \overline{\mathbb{S}}_+^{n;1}$ .  In $\mbox{O}(n^4)$ operations, Algorithm~III will either compute an optimal solution of
(\ref{eq:QP}) or determine that this problem has no finite optimal solution.  \hfill $\Box$
\end{theorem}

The extension of the above solution strategy to the case $M \in \overline{\mathbb{S}}_+^{n;k}$ for $k \geq 2$ is clear.  It suffices to notice that
for such $M$, the principal submatrix $M^{\, -i}$ belongs to $\overline{\mathbb{S}}_+^{n-1;k-1}$.  Thus Algorithm~III can be applied to these $2n$ QPs;
and the same complexity evaluation allows us to conclude that the QP (\ref{eq:QP}) can be resolved in $\mbox{O}( (2n)^k \times n^3 )$ time, for every fixed $k$.
This complexity is polynomial in $n$ when $k$ is independent of $n$.

\gap

The above solution strategy can be modified to address a related class of problems with $M \in \mathbb{S}_+^{\, n}$.
Suppose that there exists a (instead of for all) subset ${\cal I} \subset [ n ]$ (presumably of small cardinality)
such that for each $i \in {\cal I}$, the principal submatrix $M^{\, -i}$ belongs to $\overline{\mathbb{S}}_+^{n-1;k_i}$ for some (presumably small) positive integer $k_i$
and the Schur complement $( M/m_{ii} ) \in \overline{\mathbb{S}}_+^{n-1;\ell_i}$ for some (presumably small) positive integer $\ell_i$.  In this case, we can solve,
for each $i \in {\cal I}$, the two subproblems (\ref{eq:fixed lb}) and (\ref{eq:fixed ub}) defined by $M^{\, -i}$ plus an additional one defined
by $( M/m_{ii} )$.  After all these subproblems of reduced dimensions are solved,
we can either obtain an optimal solution to (\ref{eq:QP}) or conclude that the problem has no optimal solution.  If all integers $| {\cal I} |$
and $\{ k_i,\ell_i \}_{i \in {\cal I}}$ are independent of $n$, then the overall algorithm remains of polynomial complexity in terms of $n$.

\revision{\section{Some Computations}

In this section, we present numerical experiments performed to test the effectiveness of PPPA.
All experiments were conducted using {\sc matlab} R2019b on a laptop with a 2.30GHz $\text{Intel}^\text{\textregistered}$
$\text{Core}^{\text{\tiny TM}}$ i9-9880H CPU and 64 GB main memory.  The implementation used a parameter $\rho\in(0,1)$ to control
the density of the matrix $M$.  The following procedure was employed to generate instances of \eqref{eq:QP} so that $M$ is
ensured to be a member of the class $\overline{\mathbb{S}}_+^n$:

\gap

\noindent $\bullet $ Construct an $n \times n$ matrix $R$ such that each entry $R_{ij}$ ($i \geq j$) is drawn independently from
the uniform distribution $U(-0.5,0.5)$.  Let $R_{ji}=R_{ij}$ for each pair $i>j$ to ensure $M \in\mathbb{S}^n$.

\gap

\noindent $\bullet $ Construct a binary matrix $P \in \{0,1\}^{n\times n}$ such that $P_{ii}=1$ for all $i\in[n]$ and each
$P_{ij}\in\{0,1\}$ ($i> j$) is drawn independently from the Bernoulli distribution with
$\text{Prob}\left\{ P_{ij}=1\right \}=\rho$.
Let $P_{ji}=P_{ij}$ for all $i>j$.

\gap

\noindent $\bullet $ Let $M_{ij}=R_{ij}P_{ij}$ for all $i\neq j$ and $M_{ii} = R_{ii}+\sum_{j\in[n]}|M_{ij}|$ for all $i\in[n]$.
	
\gap

\noindent $\bullet $ For each $i \in [n]$, $q_i$ is drawn uniformly from [-500, 500] and the upper bound is set as
$u_i=100/\sqrt{n}$.
%
\begin{table}[H]
	\centering
	\caption{\rm Running times comparison for \pppa and \gurobi}
\begin{tabular}{cc|ccccc|c}
	\toprule\toprule
	\multirow{2}{*}{$n$}&\multirow{2}{*}{formulation}&\multicolumn{5}{c|}{density}&\\
	& &0.05 &0.1 & 0.15 &0.2&0.3&\textbf{Ave}\\\toprule
	\multirow{2}{*}{500} &       \pppa &
	0.09 &       0.11 &       0.14 &       0.14 &       0.16 &       0.13 \\
	
	&\gurobi&0.10 &       0.12 &       0.13 &       0.13 &       0.14 &       0.12 \\\midrule
	
	\multirow{2}{*}{1000} &       \pppa &0.23 &       0.33 &       0.39 &       0.44 &       0.66 &       0.41 \\
	
	&\gurobi&0.44 &       0.53 &       0.55 &       0.55 &       0.56 &       0.53 \\\midrule
	
	\multirow{2}{*}{2000} &       \pppa &0.74 &       1.21 &       1.79 &       2.49 &       5.32 &       2.31 \\
	
	&\gurobi&3.13 &       4.12 &       3.74 &       3.96 &       4.08 &       3.80 \\\midrule
%
%
%
%
%
	\multirow{2}{*}{5000} &       \pppa &5.87 &      16.82 &      35.06 &      59.84 &     117.03 &      46.92 \\
	
	&\gurobi&85.79 &      87.69 &     100.74 &      86.30 &      85.66 &      89.24 \\\midrule
	
	\multirow{2}{*}{10000} &       \pppa &52.78 &     181.24 &     377.91 &     592.19 &     989.20 &     438.66 \\
	
	&\gurobi&815.53 &     803.04 &     887.56 &     843.84 &     838.59 &     837.71 \\\midrule\bottomrule
	
\end{tabular} \label{tab:comparison}
\end{table}

Table~\ref{tab:comparison} presents the results for matrices with varying dimensions $n$ and densities $\rho$.
Five instances are generated for each combination of $n\in\{500,1000,2000,5000,10,000\}$ and $\rho\in\{0.05,0.1,0.15,0.2,0.3\}$.
For the purpose of comparison, the \gurobi method is used as a benchmark to solve all instances, which is implemented
using Gurobi~9.0 in default settings. Each row of Table~\ref{tab:comparison} reports the mean for five instances in seconds.
In addition, the statistics for different $\rho$ is averaged and shown in the last column. As one can see, \pppa outperforms \gurobi
in most problem settings except for the case of $n=500$ or $\rho=0.3$.  More importantly, as sparsity increases,
one can observe that the speedup becomes more significant and is up to almost 16 times faster (for $n =10,000$ and $\rho = 0.05$).
Particularly noteworthy is that in this case, \pppa is able to solve the problem within one minute whereas \gurobi
would take more than 12 minutes to solve it.  For the smallest $n = 500$, the times for both methods are quite small although
\gurobi is slightly faster on all the problems in this group (0.13 versus 0.12 seconds).  The experiment clearly demonstrates
the efficiency of \pppa in solving \eqref{eq:QP} as we expect in this paper.  It should be mentioned that we relied solely on
the internally available {\sc matlab} routines to carry out all the linear-algebraic calculations whereas we expect that they
are professionally implemented in \gurobi.  It is perhaps due to this reason that the solution times of PPPA start to worsen
at the last sparsity $\rho = 0.30$.  Since the main contribution of our work is not in a most efficient implementation of PPPA,
we have not attempted to improve these solution times.

\begin{figure}[H]
	\centering
	\includegraphics[width=0.7\textwidth]{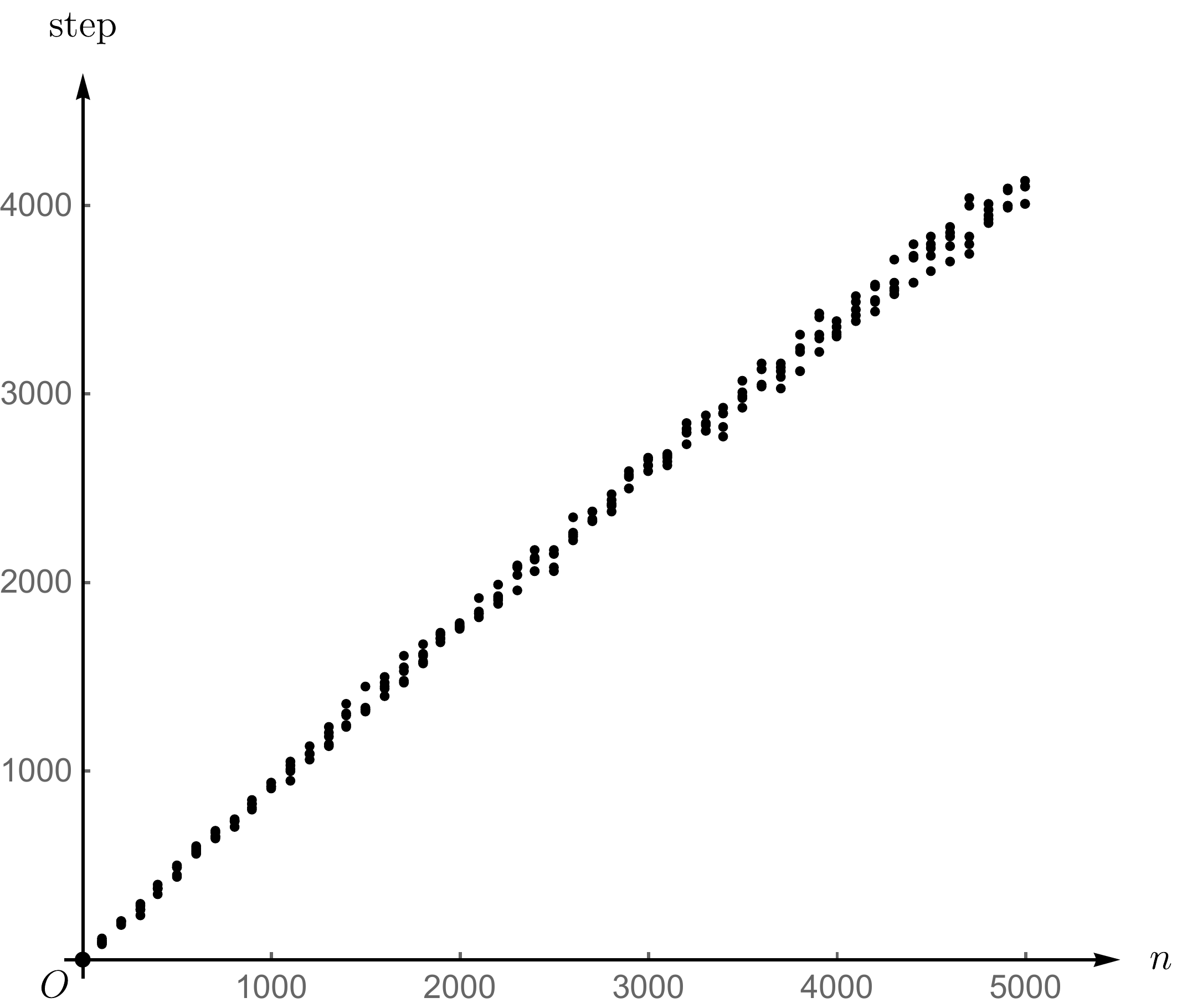}
	\caption{\rm Number of steps for \pppa with varying dimensions and fixed $\rho=0.2$}
	\label{fig:linearity}
\end{figure}

Figure~\ref{fig:linearity} shows the number of iterations in each run of \pppa as $n$ varies from  100 to 5000.
The parameter $\rho$ is fixed as $0.2$ for all instances in this experiment. As can be seen easily, the trend of the graph strongly
indicates a linear relationship between the number of steps and $n$. This matches the complexity result about \pppa as
we analyzed in Section~\ref{sec:irreducible Sbar} (see Theorem~\ref{th:irreducible barpsd}).

\section{Conclusions}

In this paper, we have analyzed the performance of the classical parametric principle pivoting algorithm on a class box-constrained quadratic programs \eqref{eq:QP}.
We prove that when the matrix belongs to $\bar{\mathbb{S}}^n$, PPPA enjoys a linear iteration complexity.
We also extend the results to the problems where the matrix is $k$-weakly quasi-diagonally dominant, and illustrate the effectiveness of PPPA with numerical experiments.
}


\begin{thebibliography}{999}

\bibitem{AdlerCottlePang16}
{\sc I.\ Adler, R.W.\ Cottle, and J.S.\ Pang}.
Some LCPs solvable in strongly polynomial time with Lemke's algorithm.
{\sl Mathematical Programming, Series A} 60(1): 477--493 (2016).

\bibitem{Chandrasekaran70}
{\sc R.\ Chandrasekaran}.
A special case of the complementary pivot problem.
{\sl Opsearch} 7: 263--268 (1970).

\bibitem{Cottle74}
{\sc R.W.\ Cottle}.
Manifestations of the Schur complement.
{\sl Linear Algebra and its Applications} 8 (1974) 189--211.

\bibitem{CottlePangStone92}
{\sc R.W.\ Cottle, J.S.\ Pang, and R.E.\ Stone}.
{\sl The Linear Complementarity Problem}.
Classics in Applied Mathematics, Volume 60. SIAM, Philadelphia (2009).
[Originally published by Academic Press, Boston (1992).]

\bibitem{GomezHePang21}
{\sc A.\ G\'omez, Z.\ He, and J.S.\ Pang}.
Linear-step solvability of some folded concave and singly-parametric sparse
optimization problems.
\revision{{\sl Mathematical Programming, Series B} (January 2022).
\url{https://doi.org/10.1007/s10107-021-01766-4}.}

\bibitem{HeHanCuiGomezPang21}
{\sc Z.\ He, S.\ Han, Y.\ Cui, A.\ G\'omez, and J.S.\ Pang}.
Comparing solution paths of sparse quadratic minimization with a Stieltjes matrix.
{\sl Optimization Online} (September 2021).

\bibitem{LiuFattahiGomezK21}
{\sc P.\ Liu, S.\ Fattahi, A.\ G\'omez, and S.\ K\"{u}\c{c}\"{u}kyavuz}.
A graph-based decomposition method for convex quadratic optimization with indicators.
\revision{{\sl Mathematical Programming} (June 2022).
\url{https://doi.org/10.1007/s10107-022-01845-0}.}

\bibitem{MonteiroAdler89}
\revision{{\sc R.D.C.\ Monteiro and I.\ Adler}.
Interior path following primal-dual algorithms. part II: Convex quadratic programming.
{\sl Mathematical Programming} 44: 43--66 (1989).
}

\bibitem{Ortega90}
{\sc J.M.\ Ortega}
{\sl Numerical Analysis: A Second Course}.
SIAM Classics in Applied Mathematics.  Volume \#3 (Philadelphia 1990).

\bibitem{PangChandrasekaran85}
{\sc J.S.\ Pang and R.\ Chandrasekaran}.
Linear complementarity problems solvable by a polynomially bounded pivoting algorithm.
{\sl Mathematical Programming Study} 25: 13--27 (1985).

\bibitem{PooleBoullion74}
{\sc G.\ Poole and T.\ Boullion}.
A survey on M-matrices.
{\sl SIAM Review} 16(4): 419--427 (1974).	

\bibitem{YeTse89}
\revision{{\sc Y.\ Ye and E.\ Tse}.
An extension of Karmarkar's projective algorithm for convex quadratic programming.
{\sl Mathematical Programming} 44: 157--179 (1989).
}
\end{thebibliography}
\end{document}